\documentclass[11pt]{article}
\usepackage{amssymb,amsthm,amsmath,amsfonts}
\usepackage[latin1]{inputenc}
\usepackage{graphicx}
\usepackage{verbatim}
\usepackage{array}
\usepackage[colorlinks=true,urlcolor=blue,
citecolor=red,linkcolor=blue,linktocpage,pdfpagelabels,
bookmarksnumbered,bookmarksopen]{hyperref}
\newtheorem{theorem}{Theorem}
\newtheorem{proposition}[theorem]{Proposition}
\newtheorem{lemma}[theorem]{Lemma}
\newtheorem{corollary}[theorem]{Corollary}

 \renewcommand{\l}{\lambda}

\renewcommand{\a}{\alpha} 
\newcommand{\g}{\gamma} \newcommand{\G}{\Gamma}
 \newcommand{\D}{\Delta}

 \renewcommand{\O}{\Omega}
\newcommand{\s}{\sigma}

\newcommand{\mf}{\mathfrak}

\newcommand{\be}{\begin{equation}}
\newcommand{\ee}{\end{equation}}

\usepackage[lmargin=2.5cm,rmargin=2.5cm,tmargin=2.5cm,bmargin=3.2cm]{geometry}

\title{ A three population Lotka-Volterra competition model with two populations interacting through an interface}
\author{Pablo \'Alvarez-Caudevilla $^1$, Cristina Br\"{a}ndle$^1$, M\'onica Molina-Becerra$^2$, Antonio Su\'arez$^3$}
\begin{document}
\pretolerance10000

\maketitle

\noindent 1. Dpto. Matem\'aticas, Universidad Carlos III de Madrid, Avda. Universidad 30, 28911, Legan\'es, Spain.

\noindent 2. Dpto. Matem\'atica Aplicada II, Escuela Polit\'ecnica Superior, Univ. de Sevilla, C/ Virgen de \'Africa, 7, 41011, Sevilla, Spain.

\noindent 3. Dpto. EDAN and IMUS. Univ. de Sevilla, Avd. Reina Mercedes, s/n, 41013,  Sevilla, Spain

\vskip 0.5cm

\noindent e-mails: pacaudev@math.uc3m.es, cbrandle@math.uc3m.es, monica@us.es, suarez@us.es.
\begin{abstract}
	In this work we consider three species competing with each other in the same habitat. One of the species lives in the entire habitat, competing with the other two species,
	while the other two inhabit two disjoint regions of the habitat. These two populations just interact on a region/interface
	which acts as a geographical barrier.  This barrier condition causes a drastic change in species behaviour compared to the classical Lotka-Volterra competitive model, showing
	very rich and new different situations depending on the several parameters involved in the system.

\end{abstract}

\noindent{\bf\small Keywords:}  {\small Interface problems, competitive systems, membrane regions interchange of
flux}

\noindent{\bf\small MSC2010:} {\small 35J60, 35J47, 35K57, 92B05}

 \section{Introduction}

\begin{figure}
	\centering
	\includegraphics[scale=0.35]{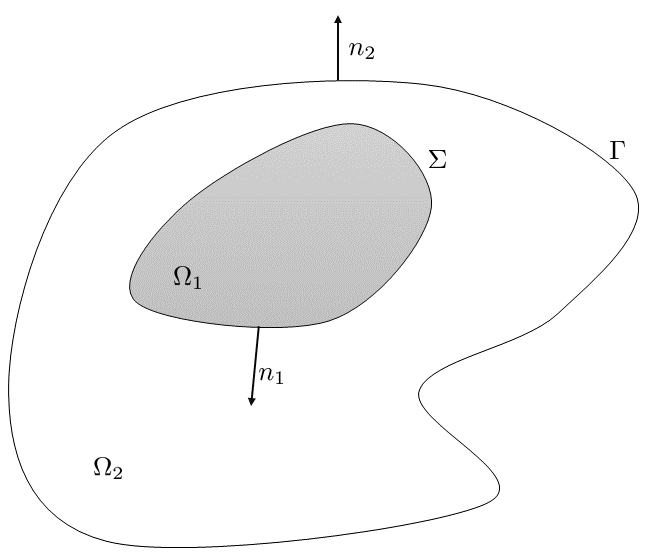}
	\caption{A possible configuration of the domain $\Omega=\Omega_1\cup \Omega_2 \cup \Sigma$.}
	\label{figure1}
\end{figure}

We aim to analyse the existence of stationary solutions for an evolution model, of different groups of populations living  in two separated regions, namely $\Omega_1$ and $\Omega_2$ of $\mathbb{R}^N$, with $N\geq 1$.
Two of them living in two different domains, whose densities are denoted by $u_1, u_2$, and have an interaction through a transition region $\Sigma$ between those two areas and the
third one $v$ living everywhere $\Omega=\Omega_1\cup\Omega_2\cup \Sigma$. In particular,
it could represent the migration/relation of two groups $u_i$ of a single population
living in $\Omega_i$, but having some effect on the other region, through an interchange of flux on the interface region $\Sigma$, that can be seen as a natural barrier, while both competing with $v$.
Thus, we analyse such a problem in the domain $\Omega=\Omega_1\cup \Omega_2 \cup \Sigma$
with $\Omega_i$ two subdomains, with internal interface $\Sigma=\partial \Omega_1$, and $\Gamma=\partial\Omega_2\setminus\Sigma$ (see Figure \ref{figure1} where we have illustrated an example of $\Omega$).
Mathematically we consider the following interaction system
\begin{equation}
\label{eq:main.system} \left\{\begin{array} {l@{\quad}l} -\Delta
u_{i} =u_{i}(\lambda_i-\alpha_i u_{i}-a_i v),&\text{ in } \Omega_i,\;\hbox{with $i=1,2$},\\ -\Delta v
=v(\mu- v-b_{1}u_{1}\chi_{\Omega_1}-b_{2}u_{2}\chi_{\Omega_2}),& \text{ in } \Omega,\\ \partial_{\bf \nu}u_{i}=\gamma_i(u_{2}-u_{1}),&\text{on $\Sigma$},\\ \partial_{\bf n} v=\partial_{\bf n}
u_{2}=0,& \text{ on $\Gamma$}. \end{array}\right.
\end{equation}
where
$a_{i}, b_{i}>0$ are positive parameters so that our model shows pure competition between
species $u_i$ and $v$. Also, $\lambda_i$ and $\mu$ are real parameters, standing for the intrinsic growth rate of $u_{i}$ and $v$ respectively, and
$ \alpha_i>0$ are referred to as the parameters representing the intrinsic competition of the species. Finally,
${\bf n_i}$ are the outward normal unitary vectors to the membrane from the domains $\Omega_i$ so that ${\bf \nu}:={\bf n_1}=-{\bf n_2}$ on $\Sigma$ and then
$$ \partial_{\bf n_1}u_2=\gamma_2(u_2-u_1) \Longleftrightarrow \partial_{\bf n_2}u_2=\gamma_2(u_1-u_2),$$
while ${\bf n}$ is the outward normal unitary vector to the domain $\Omega$.
Thus, we denote the boundary conditions by
\begin{equation}
\label{eq:BC_KK}
{\mathcal I}({\bf u})=0 \quad\mbox{on $\Sigma\cup \Gamma$}\Longleftrightarrow
\left\{\begin{array}{ll}
	\partial_{\nu}u_i=\gamma_i(u_2-u_1)& \text{on $\Sigma$},\\
	\partial_{{\bf n}}u_2=0& \text{on $\Gamma$},
\end{array}\right.
\end{equation}
where ${\bf u}=(u_1,u_2)^T$. 
We will use the subscript $1$, respectively $2$, to denote objects
(functions, parameters...) defined on $\Omega_1$ (resp. $\Omega_2$) and each time we write the subscript $i$ we mean $i=1,2$ and we will not mentioned it anymore.

%
%
%
%

The corresponding boundary conditions describing the flow through such a common boundary, the so-called membrane, are compatible with mass conservation and energy dissipation,
and called {\it Kedem-Katchalsky membrane conditions}. The conservation of mass brings to flux continuity, and the dissipation principle implies that the $L^2$-norm
of the solution is decreasing in time. This last property gives us that the flux is proportional to the difference of densities on the membrane with a proportionality coefficient,
called the membrane permeability constant, and is usually expressed as $\partial_{{\bf n}_1}
u_{1}=\gamma(u_{2}-u_{1})=\partial_{{\bf n}_1} u_{2}$, with ${\bf n}_i$ standing as the normal vectors to the membrane depending on the direction. Therefore, the main interest of system \eqref{eq:main.system}
stems from the impact of the interior boundary conditions. However, we observe that for our particular problem we assume different permeability parameters $\gamma_i$ depending on the direction of movement through the membrane.

These Kedem-Katchalsky conditions were introduced in a thermodynamic context \cite{Ked-Kat}. However, their application was also studied for Biological problems later on by
Quarteroni-Veneziani-Zunino \cite{Quarteroni} in the analysis of the dynamics of the solute in a vessel and an arterial wall. From a mathematical point of view Calabr\`o \& Zunino \cite{Cal-Zun}
obtained some results on the behaviour of a biological model for the transfer of chemicals through thin biological membranes. Moreover, in more recent works, from 2017 onwards,
there have been several articles analysing tumor invasion such as in the pressure equation in Gallinato et al.~\cite{Gallinato} or in the tumor cell density's equation in Chaplain et. al~\cite{Chaplain}.
Also, very recently by Ciavolella \& Perthame \cite{CiavolellaPaerthame} developed
a theory of weak solutions based on $L^1$ bounds assuming quadratic reaction terms.
Existence results while assuming spatial heterogeneities where obtained in \cite{PAC-Bran}. Finally,
we should mention \cite{WangSu} and \cite{Suarezetal}, in this last work the existence and uniqueness of positive solutions was analysed assuming non-symmetric boundary conditions, i.e. the permeability parameter $\gamma$
depends on the direction of the flux.

\subsection{Main results}

For problem \eqref{eq:main.system} we analyse the existence of solutions in terms of the parameters involved in the system, i.e. $\lambda_i$ and $\mu$. This new approach opens up a wide range of possibilities for such an existence of
solutions. Also, providing us with new unexpected behaviours, some of them clearly due to the effect of the membrane interface region.

It is clear what happens just by analysing the logistic problem, just assuming the components $u_1$ and $u_2$ in  \eqref{eq:main.system}. In other words, somehow assuming semitrivial solutions with $v=0$, under Kedem-Katchalsky conditions
\eqref{eq:BC_KK}. So that we find the existence of semitrivial positive solutions of the form $(u_1,u_1,0)$ if and only if $\Lambda_1(-\lambda_1,-\lambda_2)<0$, where $\Lambda_1(-\lambda_1,-\lambda_2)$ stands for
the principal eigenvalue of the interface problem
$$
	\left\{\begin{array} {l@{\quad}l}
		-\Delta \varphi_i - \lambda_i
		\varphi_i=\nu \varphi_i&  \hbox{in}\quad \Omega_i,\\
		\mathcal{I}({\bf \Phi})=0 &\hbox{on } \Sigma\cup\Gamma,
		\end{array}\right.
$$
with ${\bf \Phi}=(\varphi_1,\varphi_2)^T$.
Moreover, we also have semitrivial solutions of the form $({\bf 0},v)^T$ for a logistic type equation of the form
\begin{equation} \label{eq:logisticgen}
\left\{\begin{array} {l@{\quad}l} -\Delta v +c(x)v=\mu v- \beta v^2&
\mbox{in $\Omega$},\\
\mathcal{B} v=0&\mbox{on $\partial\Omega$,}
\end{array}\right.
\end{equation}
where $\mu\in\mathbb{R}$, $\beta>0$
and $c(x)\in L^\infty(\Omega)$, for general homogeneous boundary conditions $\mathcal{B} v=0$ on $\partial\Omega$, if $\mu >\sigma_1^{\Omega}[-\Delta+c(x);\mathcal{B}]$, where $\sigma_1^{\Omega}[-\Delta+c(x);\mathcal{B}]$ denotes the principal eigenvalue of the problem
$$
\left\{\begin{array} {l@{\quad}l} -\Delta v +c(x)v=\lambda v&
\mbox{in $\Omega$},\\
\mathcal{B} v=0&\mbox{on $\partial\Omega$.}
\end{array}\right.
$$
Observe that, in case $\mathcal{B} v=0$ stands for homogeneous Neumann boundary conditions, that eigenvalue lower bound for the existence of positive solutions would be 0, hence, we have that $\mu>0$ for the
existence of these types of semitrivial solutions.

Thus, here we also show that the  necessary condition that guarantees the existence of {\it coexistence states} $(u_1,u_2,v)^T$ of problem \eqref{eq:main.system}  (that is, solutions for which $u_i$ and $v$ are positive) is
$$
\Lambda_1(-\lambda_1,-\lambda_2)<0\quad \hbox{and}\quad \mu>0.
$$
From those results we show that once we have fixed the values of the parameters $\lambda_1$ and $\lambda_2$,  such that $\Lambda_1(-\lambda_1,-\lambda_2)<0$, i.e. the existence of semitrivial solutions
given by $(\theta_{\lambda_1},\theta_{\lambda_2},0)^T$ is guaranteed, there exits a value for the parameter $\mu$, denoted in this paper by $\mu_*(\lambda_1,\lambda_2)$, such that when $\mu>\mu_*(\lambda_1,\lambda_2)$,
problem \eqref{eq:main.system} does not possess a coexistence state.

Furthermore, for those fixed values of the parameters  $\lambda_1$ and $\lambda_2$, such that $\Lambda_1(-\lambda_1,-\lambda_2)<0$, we additionally prove in this paper the existence of a bifurcation point $\mu_0$ such
that a continuum of coexistence states emanates from the semitrivial solution $(\theta_{\lambda_1},\theta_{\lambda_2},0)^T$, up to another value $\mu_1$. Hence, having the existence of coexistence states when
$\mu$ is between those two values $\mu_0$ and $\mu_1$. As shown below in the analysis performed in this paper, such a continuum of coexistence states, denoted by $\mathcal{C}$ satisfies that
$\mbox{Proj}_\mathbb{R}(\mathcal{C})\subset[0,\mu_*(\lambda_1,\lambda_2)]$, so that $\mu_0,\mu_1\in (0,\mu_*(\lambda_1,\lambda_2))$.

Subsequently, to study the region of coexistence states in terms of the main parameters involved in the system, let us say $\lambda_i$ and $\mu$, we shall construct several auxiliary functions providing us with
equivalent conditions for such existence of coexistence states.

In particular, assuming without loss of generality that $\mu_0<\mu_1$ we establish an equivalence between the coexistence states given by $\mu\in (\mu_0,\mu_1)$, for $\lambda_1$ and $\lambda_2$ fixed, with $\Lambda_1(-\lambda_1,-\lambda_2)<0$,
as discussed above, with the expression
\begin{equation}
	\label{condicion}
	(\mu-G(\lambda_1,\lambda_2))\cdot \Lambda_1(-\lambda_1+a_1\mu,-\lambda_2+a_2\mu)<0,
\end{equation}
where the function $G(\lambda_1,\lambda_2)$ is denoted by
\begin{equation}
\label{eq:G_lambdas}
G(\lambda_1,\lambda_2):= \sigma_1^\Omega[-\Delta + b_1 \theta_{\lambda_1}\chi_{\Omega_1}+  b_2 \theta_{\lambda_2} \chi_{\Omega_2};\mathcal{N}],
\end{equation}
where $\mathcal{N}$ denotes homogeneous Neumann boundary condition. In other words, we will have coexistence states if
$$\left\{\begin{array}{l} \hbox{either}, \quad \mu-G(\lambda_1,\lambda_2)<0\quad \hbox{and}\quad \Lambda_1(-\lambda_1+a_1\mu,-\lambda_2+a_2\mu)>0,\\
\hbox{or}, \quad \mu-G(\lambda_1,\lambda_2)>0\quad \hbox{and}\quad \Lambda_1(-\lambda_1+a_1\mu,-\lambda_2+a_2\mu)<0.\end{array}\right.$$
First we observe that the value of the parameter $\mu$ from where a continuum $\mathcal{C}$ of coexistence states emanates from the semitrivial solutions  $(\theta_{\lambda_1},\theta_{\lambda_2},0)^T$ is actually
$$\mu_0=\sigma_1^\Omega[-\Delta + b_1 \theta_{\lambda_1}\chi_{\Omega_1}+  b_2 \theta_{\lambda_2} \chi_{\Omega_2};\mathcal{N}].$$
Hence, the analysis of $\mu-G(\lambda_1,\lambda_2)$ is equivalent to study $\mu-\mu_0(\lambda_1,\lambda_2)$ for different values of the parameters $\lambda_i$.
Thus, to analyse such an expression our approach is based on fixing the values of one of the parameters, in the sequel we fix $\lambda_2$, and analyse a function denoted by $g(\lambda_1)\equiv G(\lambda_1,\lambda_2)$ once $\lambda_2$ is fixed.

In doing so, we study the behaviour at infinity of the function $g(\lambda_1)$. Proving then, that such a function is continuous and increasing in terms of $\lambda_1$ and bounded above when $\lambda_1\to +\infty$. On the other hand, the
behaviour when $\lambda_1$ goes to $-\infty$ the limit will depend on the election of the fixed value $\lambda_2$. Having a lower bound if $\lambda_2 \geq \sigma_2$, for a particular value $\sigma_2$ shown below. While, if $\lambda_2<\sigma_2$
the function $g(\lambda_1)$ is only defined for values of $\lambda_1$ sufficiently large.

At this point it is crucial to mention the very important result we have proved in Lemma\;\ref{autovalor} which is key in obtaining the existence of coexistence states. In that lemma an auxiliary function
$\mathcal{H}:(-\infty,\sigma_2)\mapsto \mathbb{R}$ is studied.  
Consequently, from Lemma\;\ref{autovalor} we are able to establish that, for such a fixed $\lambda_2$:
$$\hbox{if $\lambda_2\geq \sigma_2 \Rightarrow \Lambda_1(-\lambda_1,-\lambda_2)<0$, for any $\lambda_1\in \mathbb{R}$},$$
on the other hand,
\begin{equation}
	\label{eq:con_lambda2}
	\hbox{if $\lambda_2< \sigma_2 \Rightarrow \Lambda_1(-\lambda_1,-\lambda_2)<0$, if $\lambda_1> \mathcal{H}(\lambda_2)$}.
	\end{equation}
Furthermore, analysing now the sign of the eigenvalue $\Lambda_1(-\lambda_1+a_1\mu,-\lambda_2+a_2\mu)$ we find an equivalence for this  function, just applying Lemma\;\ref{autovalor}, having an important value
$\sigma_1>-\lambda_2+a_2\mu$ which will distinguish the behaviour of this eigenvalue. Indeed, if $\mu>(\lambda_2-\sigma_2)/a_2$ we find that
$$
\Lambda_1(-\l_1+a_1\mu,-\l_2+a_2\mu)<0\quad\mbox{if and only if}\quad \l_1>\mathcal{G}(\mu),
$$
and
\begin{equation}
	\label{eq:con_Lambda1}
\Lambda_1(-\l_1+a_1\mu,-\l_2+a_2\mu)>0\quad\mbox{if and only if}\quad \lambda_1<\mathcal{G}(\mu),
\end{equation}
having an equivalence between the sign of the eigenvalue $\Lambda_1(-\lambda_1+a_1\mu,-\lambda_2+a_2\mu)$ and the map $\mathcal{G}(\mu)=a_1\mu+\mathcal{H}(\lambda_2-a_2\mu)$ defined for values of $\mu$
larger than $(\lambda_2-\sigma_2)/a_2$.

Therefore, we can now present our results, that characterize the region of coexistence states $\mathcal{R}$ of the system, in particular, fixing $\lambda_2$ and moving $\lambda_1$ and $\mu$. Thus, and thanks to what we have discussed above,
condition \eqref{condicion}, once we fix $\lambda_2$, is actually equivalent to
\begin{equation}
	\label{condicion2}
	(\mu-g(\l_1))\cdot(\l_1-\mathcal{G}(\mu))<0,
\end{equation}
for system (\ref{eq:main.system}) possessing at least a coexistence state. 

First of all, let us note that to study the coexistence region $\mathcal{R}$ we have to ensure that $\mu>0$ and $\Lambda_1(-\lambda_1,-\lambda_2)< 0$. Otherwise there will not be coexistence states, as discussed above. Then $\mathcal{R}$ will be determined by the curves given in~\eqref{condicion2}. The following discussion can be seen thanks to Figures~\ref{fig1} and~\ref{fig2} in Section~\ref{sect:coexistence}, where we have graphed the functions $\mu=g(\lambda_1)$ and $\lambda_1=\mathcal{G}(\mu)$ in different situations.

Fixing $\lambda_2$ so that $\lambda_2\geq \sigma_2>0$,  see Figure~\ref{fig2}, we find the existence of coexistence states for positive values of $\mu$,  and, surprisingly, for  negative values of $\lambda_1$, since $\Lambda_1(-\lambda_1,-\lambda_2)<0$ for all $\lambda_1$. Moreover, the function $g$ is defined for all the vales of $\lambda_1$ and $\mathcal{R}$ will be determined by the asymptotic behaviour of $g$.
 We observe that this situation is not classical at all and it is basically due to
the effect of the interface. Indeed, since thanks to the interface condition there is population coming through the interface $\Sigma$,  even when the parameter $\lambda_1$ is negative, very negative indeed, the species still coexist. Therefore, we have found
 here a situation
where the interface produces a cooperative effect despite having negative growth rates for some of the populations. In that sense the existence of solutions is recovered for negative $\lambda_1$ while $\lambda_2$ and $\mu$ are positive.

On the other hand, if we now fix $\lambda_2$ such that $\lambda_2<\sigma_2$, see Figure~\ref{fig1}, assuming that the parameter $\mu$ is positive, we will find coexistence states for positive large values of the parameter $\lambda_1$,
in particular, if $\lambda_1 \in [\mathcal{H}(\lambda_2),\infty)$, even
having negative values of the parameter $\lambda_2$. 
Moreover, we prove in this paper that choosing $\lambda_2$ negative
is equivalent to have that $\mathcal{H}(\lambda_2)>0$.
Also, if $\lambda_2<\sigma_2$ but positive, so that $\mathcal{H}(\lambda_2)<0$, we find coexistence states for negative values of the parameter $\lambda_1$ as well. In summary, we again have  the existence of coexistence states when one of the parameters
is positive while the other is negative.

Up to what we have explained so far, from the results obtained in this paper we can assure that (see Figures \ref{fig1} and \ref{fig2}):
\begin{enumerate}
	\item If $\mu\leq 0$, then (\ref{eq:main.system}) does not possess coexistence states. 
	\item There exists a positive value $g^\infty$, that depends on $\l_2$, such that if $\mu>g^\infty$ then  (\ref{eq:main.system}) possess  at least a coexistence state if   $\lambda_1>\mathcal{G}(\mu)$.
	\item If $0<\mu\leq g^\infty$ and $\l_2<\sigma_2$, there exists a function $\mathcal{H}$ such that for $\lambda_1<\mathcal{H}(\lambda_2)$, system~\eqref{eq:main.system} does not possess coexistence states. Moreover,   (\ref{eq:main.system}) possess  at least a coexistence state if   $\lambda_1\in ( \underline{\l}_1,\overline{\l}_1)$ with $\underline{\l}_1>0$ if $\l_2<0$ and  $\underline{\l}_1<0$ if $\l_2>0$. In this last case, the competitive species coexist  even for negative values of the growth rate of $u_1$.
	\item If $0<\mu\leq g^\infty$ and $\l_2\geq \sigma_2$, then there exists $g_\infty$, depending on $\l_2$, such that for 
	$$\mu\in(\min\{(\l_2-\s_2)/a_2,g_\infty\}, \max\{(\l_2-\s_2)/a_2,g_\infty\})
	$$
	then (\ref{eq:main.system}) possess  at least a coexistence state if   $\lambda_1<\underline{\l}_{1,*}$ for some value $\underline{\l}_{1,*}$. Again, the species can coexist for negative values of the growth rate of the species $u_1$. 
		 

	
\item Principle of Competitive Exclusion: Fixed $\lambda_1$ and $\lambda_2$ so that $\Lambda_1(-\lambda_1,-\lambda_2)<0$, then~(\ref{eq:main.system}) does not possess coexistence states for $\mu$ large, $\mu\geq \mu_*(\lambda_1,\lambda_2)$.
	
		\item Fix $\mu\geq \s_2$, then (\ref{eq:main.system})  possesses at least a coexistence state even for $\l_1$ large and positive, see Section~\ref{sect:limit}. Hence, the Principle of Competitive Exclusion is broken in this situation.
		
\end{enumerate}
In conclusion, the membrane assumed for the competition system (\ref{eq:main.system}) provides us with a richer existence of solutions from the classical competition problem where there is not such an interface membrane region. Indeed, we actually find that the Principle of Competitive Exclusion is broken for certain cases when the parameter $\lambda_1$ and $\mu$ are sufficiently large. In addition, we obtain as well several other situations where system
(\ref{eq:main.system}) possesses at least a coexistence state for negative values of the growth rate of the species $u_1$.   This is a completely new situation compared to the case of competition, where the growth coefficients of each species must be positive for the species to coexist. 

Finally, we also analyse the case when the parameters $\lambda_1$ and $\lambda_2$ are assumed to be the same, i.e. $\lambda_1=\lambda_2=\lambda$. In this particular case we observe that,
following a similar analysis as the one previously performed for different parameters, the system shows a more classical Lotka-Volterra system, with the interface $\Sigma$ not showing such an important effect as we have observed for the previous one (see Figure \ref{fig3}).
Such a conclusion makes sense, since in this particular case, we cannot have different signs for the parameters $\lambda_1$ and $\lambda_2$ and $u_1$ and $u_2$ behave in some sense as if they where a unique population.

\section{Preliminaries}
In this section we fix notation concerning eigenvalue problems, that we will use through the paper. Properties of these eigenvalues as well as existence and uniqueness results of the associated problems are studied in~\cite{AlvarezBrandleMolinaSuarez}. Moreover, we state the properties of semitrivial solutions of~\eqref{eq:main.system}.

\subsection{Eigenvalue problems}
Let us consider a general domain $D\subset \mathbb{R}^N$ which is regular and bounded, such that
$$
\partial D=\Gamma_1\cup \Gamma_2,
$$
where $\Gamma_1$ and $\Gamma_2$ are two disjoint open and closed subsets of the boundary.
Note that, it is possible that some $\Gamma_i=\emptyset$ for some $i=1,2$.
Now, consider the eigenvalue problem
\begin{equation}
\label{eq:eign_prob_logistic}
 \left\{\begin{array} {l@{\quad}l} -d\Delta \varphi + c(x)
\varphi=\rho \varphi,& \hbox{in} \quad D,\\ \mathcal{B} \varphi=0,& \hbox{on}\quad \partial D.
\end{array}\right.
\end{equation}
with $d>0$, $c\in C(\overline D)$ and
\begin{equation}
\label{eq:BC_eign_prob}
{\mathcal{B}}\varphi:=
{\partial_{\bf n_i}} \varphi +g_i \varphi \quad  \mbox{on $\Gamma_i$,}
\end{equation}
where ${\bf n_i}$ is a nowhere tangent vector field and $g_i\in C(\Gamma_i)$.  We denote the principal eigenvalue of~\eqref{eq:eign_prob_logistic}  as
\begin{equation}
\label{eq:princ_eign}
\sigma_1^D[-d\Delta+c;{\mathcal{B}}]=\sigma_1^D[-d\Delta+c;\mathcal{N}+g_1;\mathcal{N}+g_2]
\end{equation}
with $\mathcal{N}$ standing for the normal derivative shown in the general boundary conditions \eqref{eq:BC_eign_prob}. Moreover, for completion in the sequel, we will use $\mathcal{D}$ to denote Dirichlet boundary data.
As mentioned, properties $\sigma_1^D[-d\Delta+c;{\mathcal{B}}]$ are stated in~\cite{Santi}, see also ~\cite{AlvarezBrandleMolinaSuarez}. In particular, there it is shown that $\sigma_1^D[-d\Delta+c;{\mathcal{B}}]$ is continuous and increasing in $c$ and the limits when $d$ tends to $0$ or $\infty$ are characterized.
%

Consider also the interface eigenvalue system given by
\begin{equation}
	\label{eq:eign_prob_logistic.system}
	\left\{\begin{array} {l@{\quad}l}
		-d\Delta \varphi_i + c_i(x)
		\varphi_i=\nu \varphi_i&  \hbox{in}\quad \Omega_i,\\
		\mathcal{I}({\bf \Phi})=0 &\hbox{on } \Sigma\cup\Gamma,
		\end{array}\right.
\end{equation}
where we assume that  $c_i\in C(\overline\Omega_i)$. We denote by
$$\Lambda_1(-d\Delta+c_1,-d\Delta+c_2),$$ the principal eigenvalue of problem (\ref{eq:eign_prob_logistic.system}) or just
$$
\Lambda_1(c_1,c_2):=\Lambda_1(-\Delta+c_1,-\Delta+c_2),
$$
 when $d=1$.

 As for the scalar problem, in~\cite{AlvarezBrandleMolinaSuarez} the continuity and monotonicity of $\Lambda_1$ is proved. Moreover, the following limit behaviours can be deduced, that will be crucial in our sequel analysis, are also stated:
  \begin{equation}
\label{dainfinito}
\lim_{d\to \infty} \Lambda_1(-d\Delta+c_{1,d},-d\Delta+c_{2,d})= \frac{\displaystyle\gamma_2\int_{\Omega_1}c_1+\gamma_1\int_{\Omega_2}c_2}{\gamma_2|\Omega_1|+\gamma_1|\Omega_2|},
\quad\text{where $c_{i,d}\to c_i$ in $L^\infty(\Omega_i)$}
\end{equation}
and
  \begin{equation}
\label{prop.lambda.inf}
\lim_{\mu\to\infty} \Lambda_1(-d\D+\mu c_1,-d\D+\mu c_2)=\infty,\quad\text{for $c_i>0$ in $\Omega_i$}.
\end{equation}
\subsection{Semi-trivial solutions}
It is clear that the trivial solution $(0,0,0)$ verifies system \eqref{eq:main.system}. Focusing on the characterization of semitrivial solutions we will consider them to be either of the form $(0,0,v)$ or $(u_1,u_2,0)$.

In the absence of the species $u_1$ and $u_2$,
the population $v$ in~\eqref{eq:main.system} behaves following a logistic type equation under homogeneous Neumann boundary conditions of the form:
\begin{equation} \label{eq:logisticpar}
\left\{\begin{array} {l@{\quad}l} -\Delta v =\mu v-  v^2&
\mbox{in $\Omega$},\\
\partial_{\bf n} v=0&\mbox{on $\partial\Omega$,}
\end{array}\right.
\end{equation}
where $\mu\in\mathbb{R}$. A more general version of that equation, assuming Robin boundary conditions, is \eqref{eq:logisticgen}.
 Thus, we find the following result.
\begin{theorem}
\label{Theo_logistic_classic}
There exists a unique  positive solution
of~\eqref{eq:logisticgen}, denoted by $\eta_{\mu,c}$, if and only if $\mu >\sigma_1^{\Omega}[-\Delta
+c(x);\mathcal{B}]$. Moreover,
$$
\frac{\mu-\sigma_1^\Omega[-\Delta +c(x);\mathcal{B}]}{\beta\|\varphi_1\|_\infty}\varphi_1(x)\leq \eta_{\mu,c}(x)\leq \frac{\mu-c_L}{\beta},
$$
where $\varphi_1$ is a positive eigenfunction associated to $\sigma_1^\Omega[-\Delta +c(x);\mathcal{B}]$ and $c_L=\min_{x\in \overline\Omega} c(x)$.
\end{theorem}

The proof of Theorem\;\ref{Theo_logistic_classic} has been performed by several authors during the last decades and it is a very well known result, see for instance \cite{Santi2}.

In particular, we point out, for further purposes, that we introduce the notation $\eta_{\mu}$ as the solution
of~\eqref{eq:logisticgen} with $c\equiv
0$, $\beta=1$ and $\mathcal{B}v=\mathcal{N}v=\partial_{{\bf n}} v$. That is, $\eta_{\mu}$ is the solution of~\eqref{eq:logisticpar} and
$
\eta_\mu=\mu
$.

On the other hand, when $v=0$, semitrivial solutions $(u_1,u_2,0)$ of~\eqref{eq:main.system} are solutions to
\begin{equation}
\label{eq:membrane.logis} \left\{\begin{array} {l@{\quad}l} -\Delta
u_{i} =u_{i}(\lambda_i-\alpha_i u_{i}),& \hbox{in}\quad \Omega_i,\\ \mathcal{I}({\bf u})=0, & \hbox{on}\quad \Sigma\cup \Gamma.
\end{array}\right.
\end{equation}
This system has been considered in~\cite{AlvarezBrandleMolinaSuarez} both for $\lambda_1=\lambda_2$ and $\lambda_1\neq \lambda_2$.

Let us denote by $\omega_1$ and $\omega_2$ the respective solutions of the uncoupled problems
\begin{equation}
\label{eq:logisuna}
\left\{\begin{array} {l@{\quad}l} -\Delta
u_{1} =u_{1}(\lambda_1-\alpha_1 u_{1}),& \hbox{in}\quad \Omega_1,\\
 \partial_{\bf n_1}
u_{1}+ \gamma _1u_{1}=0,&  \hbox{on}\quad \Sigma, \end{array}\right.
\end{equation} and
\begin{equation}
\label{omega2}
\left\{
\begin{array}{ll}
 -\D u_2=u_2(\l_2-\alpha_2u_2) & \mbox{in $\O_2$,}\\
 \partial_{\bf n_2} u_2+\g_2u_2=0 &  \mbox{on $\Sigma$,}\\
 \partial_{\bf n} u_2=0 &  \mbox{on $\Gamma$.}\\
\end{array}
\right.
\end{equation}
Also, we define the principal eigenvalues for the operators shown in problems \eqref{eq:logisuna} and \eqref{omega2} above, under their respective boundary conditions as
\begin{equation}
  \label{eq:def.sigmas}
\sigma_1:= \sigma_1^{\Omega_1}[-\D;\mathcal{N}+\gamma_1],\qquad \sigma_2:= \sigma_1^{\Omega_2}[-\D;\mathcal{N}+\gamma_2;\mathcal{N}].
\end{equation}
Observe that $\sigma_1,\sigma_2>0$. Thanks to Theorem \ref{Theo_logistic_classic}, we conclude that each $\omega_i$ exists and are positive  if and only if
$\l_i>\s_i$. For the semitrivial solution $(u_1,u_2,0)$ of~\eqref{eq:main.system} we show the next result, which provides us with a condition for the existence of positive solutions
in terms of the principal eigenvalue of problem \eqref{eq:eign_prob_logistic.system} and bounds for the solutions as well.
\begin{theorem}
There exists a unique  positive solution
of~{\rm\eqref{eq:membrane.logis}}, denoted by $(\theta_{\l_1},\theta_{\l_2})$, if and only if
\begin{equation}
\label{condinecesufi}
\Lambda_1(-\lambda_1,-\lambda_2)<0.
\end{equation}
Moreover,
\begin{equation}
\label{cotaimpor}
\omega_i\leq \theta_{\l_i}\leq \max\left\{\frac{\l_1}{(\alpha_1)_L},\frac{\l_2}{(\alpha_2)_L}\right\}\quad\mbox{in $\O_i$,}
\end{equation}
where $(\alpha_1)_L=\min_{x\in \overline\Omega_i} \alpha_i(x)$.
\end{theorem}

\begin{proof}
  See~\cite{AlvarezBrandleMolinaSuarez} for a complete proof.
\end{proof}

\begin{proposition}
\label{exislogismem}
The map $(\l_1,\l_2)\mapsto (\theta_{\l_1},\theta_{\l_2})\in L^\infty(\O_1)\times L^\infty(\O_2)$ is increasing and continuous. Moreover,
\begin{equation}
\label{eigenposi}
\Lambda_1(-\l_1+2\theta_{\l_1}\a_1,-\l_2+2\theta_{\l_2}\a_2)>0.
\end{equation}
\end{proposition}
\begin{proof}
Take $\l_1\leq \overline\l_1$ and $\l_2\leq \overline\l_2$. Then, it is clear that $(\underline u_1,\underline u_2)=(\theta_{\l_1},\theta_{\l_2})$ is a subsolution of~\eqref{eq:membrane.logis} with $\l_1=\overline\l_1$ and $\l_2=\overline\l_2$ and
$(\overline u_1,\overline u_2)=(K,K)$ is a supersolution, where $K$ is a sufficiently large positive constant. Then,
$$
\theta_{\l_1}\leq \theta_{\overline\l_1} \quad\mbox{in $\O_1$ \quad and}\quad\theta_{\l_2}\leq \theta_{\overline\l_2}\quad\mbox{in $\O_2$.}
$$
The continuity follows by an standard argument and the elliptic regularity.

Finally, using that
$$
0=\Lambda_1(-\l_1+\theta_{\l_1}\a_1,-\l_2+\theta_{\l_2}\a_2)
$$
and that $\Lambda_1$ is increasing in $(c_1, c_2)$, see~\cite[Lemma 3.1]{AlvarezBrandleMolinaSuarez}, we obtain~\eqref{eigenposi}.
\end{proof}
We also include the following result which will be applied in Section 4 and whose proof is shown in~\cite{AlvarezBrandleMolinaSuarez}. To this aim, let us define $L_{\lambda_2}$ the so-called large solution of the problem
\begin{equation}
  \label{eq:large}
 \left\{\begin{array} {l@{\quad}l} -\Delta u =\lambda_2 u-\alpha_2 u^2,& \hbox{in}\quad \Omega_2,\\
 u=\infty,& \hbox{on}\quad \Sigma,\\
 \partial_{\bf n} u=0, & \hbox{on}\quad \Gamma,
\end{array}\right.
\end{equation}
where the boundary data on $\Sigma$ has to be understood as $u(x)\to \infty$ as ${\rm dist}(x,\Gamma)\to 0$. The existence of this kind of solution is shown in~\cite{DuHuang} and~\cite{LG-book}.%

\begin{proposition}\label{Prop_v0}
Consider a fixed value of the parameter $\lambda_2$ such that $(\theta_{\l_1},\theta_{\l_2})$ is a positive solution of problem \eqref{eq:membrane.logis}. Then,
\begin{equation}
\label{limit:lambda_2}
\lim_{\lambda_1\to \infty} \theta_{\l_1}(x)=+\infty,\quad \hbox{for all $x\in \overline \Omega_1$}, \quad\text{and}\quad \lim_{\lambda_1\to \infty} \theta_{\l_2}=L_{\lambda_2},\quad \hbox{in $C^2(\Omega_2).$}
\end{equation}
\end{proposition}





\section{Coexistence states}
In this section we characterize the existence and non-existence of positive solutions (coexistence states) of system~\eqref{eq:main.system} in terms of $\mu$, when $\lambda_1$ and $\lambda_2$ are fixed.

\subsection{A priori bounds}
First of all we prove a priori bounds of a coexistence state of problem~\eqref{eq:main.system}.

\begin{theorem}
\label{cotassolu}
Let $(u_1,u_2, v)$ be a coexistence state  of \eqref{eq:main.system}. Then,
\begin{equation}
\label{cotas}
\|u_i\|_{C(\overline\Omega_i)}\leq \max\left\{\frac{\l_1}{(\alpha_1)_L},\frac{\l_2}{(\alpha_2)_L}\right\}, \quad \hbox{and} \quad \|v\|_{C(\overline\Omega)} \leq \mu.
\end{equation}
\end{theorem}
\begin{proof}
Let $(u_1,u_2, v)$ a coexistence state  of \eqref{eq:main.system}, then, by construction and since we are assuming positive solutions, $(u_1,u_2)$ is a sub-solution of \eqref{eq:membrane.logis}. Hence, we conclude that
$$
u_i\leq \theta_{\l_i}\leq \max\left\{\frac{\l_1}{(\alpha_1)_L},\frac{\l_2}{(\alpha_2)_L}\right\},
$$
where we have used~\eqref{cotaimpor} for the second inequality.
On the other hand, there exists
$x_M$ in $\Omega$ such that
$$
 v(x_M)=\|v\|_{C(\overline\Omega)}.$$
Thus, using Lemma 2.1 in \cite{Lou-ni} we obtain that
$$v(x_M)(\mu-v(x_M)- b_{1}u_{1}(x_M)\chi_{\Omega_1} -b_{2}u_{2}(x_M)\chi_{\Omega_2})\geq 0,$$
which turns out to be true if  $\mu\geq v(x_M)+ b_{1}u_{1}(x_M)\chi_{\Omega_1} +b_{2}u_{2}(x_M)\chi_{\Omega_2}$ and  hence
$
 \|v\|_{C(\overline\Omega)}\leq \mu.$
\end{proof}

\subsection{Existence and non-existence of coexistence states}

As a direct consequence of the previous result we get a necessary condition for the existence of coexistence states.
\begin{corollary}
\label{col:nece_cond}
  If there exists a coexistence state $(u_1,u_2,v)$ of  \eqref{eq:main.system}, then
\begin{equation*}
\label{con_exis_lambda}
\Lambda_1(-\lambda_1,-\lambda_2)<0, \quad \hbox{and} \quad \mu>0.
\end{equation*}
\end{corollary}
Furthermore, the next result shows conditions under which \eqref{eq:main.system} does not have coexistence states.
\begin{proposition}
\label{nomularge}
Let $\lambda_1$ and $\lambda_2$ fixed such that  $\Lambda_1(-\lambda_1,-\lambda_2)<0$. Then,
there exists $\mu_*=\mu_*(\l_1,\l_2)$ such that \eqref{eq:main.system} does not have coexistence states  if $\mu>\mu_*$.
\end{proposition}
\begin{proof}
Assume by contradiction that there exists a coexistence state $(u_1,u_2, v)$ of \eqref{eq:main.system} for $\mu$ large enough.
Without loss of generality, we can assume that $\l_1/(\alpha_1)_{L}\leq \l_2/(\alpha_2)_{L}$. Then, by~\eqref{cotas} we have that $u_i\leq \l_2/(\alpha_2)_{L}$ and
$$
-\Delta v\geq v(\mu -v-b_{1}( \lambda_2/(\alpha_2)_{L})\chi_{\Omega_1}-b_{2}( \lambda_2/(\alpha_2)_{L})\chi_{\Omega_2}))\quad\mbox{in $\O$.}
$$
Hence, by comparison,
$$
v\geq \omega_{\mu,c}
$$
where $\omega_{\mu,c}$ is the unique solution of (\ref{eq:logisticgen}) with
$$c(x)=( \l_2/(\alpha_2)_{L})(b_{1}\chi_{\Omega_1}+b_{2}\chi_{\Omega_2}),$$ which exists for $\mu$ large enough, see Theorem~\ref{Theo_logistic_classic}. Consequently,
$$
0=\Lambda_1(-\l_1+\alpha_1u_1+a_1v,-\l_2+\alpha_2u_2+a_2v)>\Lambda_1(-\l_1+a_1\omega_{\mu,c},-\l_2+a_2\omega_{\mu,c}u),
$$
and this is a contradiction for $\mu$ large. Indeed, since due to Theorem \ref{Theo_logistic_classic} we get that
$$
\frac{\mu-\sigma_1^\Omega[-\Delta-c;\mathcal{N}]}{\|\varphi_1\|_\infty}\varphi_1(x)\leq \omega_{\mu,c},
$$
we  conclude, see~\eqref{prop.lambda.inf}, that
$$
\Lambda_1(-\l_1+a_1\omega_\mu,-\l_2+a_2\omega_\mu)\to \infty\quad\mbox{as $\mu\to\infty,$}
$$
arriving at a contradiction.
\end{proof}

Next, we show the existence of coexistence states depending on the value of the  parameters. We fix the values of $\lambda_i$ so that the existence of a semitrivial solution of the form $(\theta_{\lambda_1},\theta_{\lambda_2},0)$ is guaranteed. This solution depends on $\mu$ and  we denote it by ${\bf S_\mu}:=(\mu,(\theta_{\lambda_1},\theta_{\lambda_2},0))$.

\begin{theorem}
\label{Theo_Bif}
Let  $\lambda_i$ be fixed such that $\Lambda_1(-\l_1,-\l_2)<0$ and ${\bf S_\mu}$ the associated semitrivial solution of~\eqref{eq:main.system}. Then there exists a bifurcation point $\mu_0=\mu_0(\l_1,\l_2)$, given by
\begin{equation}
\label{con:mu_zero}
\mu_0=\sigma_1^\Omega[-\Delta  +b_{1}  \theta_{\lambda_1} \chi_{\Omega_1}+b_{2}  \theta_{\lambda_2}\chi_{\Omega_2};\mathcal{N}],
\end{equation}
 such that a smooth curve of coexistence states emanates from the semitrivial solution ${\bf S_\mu}$.
\end{theorem}



\begin{proof}
Take $\rho>0$, so that the corresponding eigenvalue of the problem \eqref{eq:eign_prob_logistic.system} satisfies $\Lambda_1(\rho,\rho)>0$. Under this assumption, the linear system
$$
-\Delta u_i+\rho u_i=f_i(x)\quad\mbox{in $\O_i$,}\quad \mathcal{I}({\bf u})=0\quad\mbox{on $\G\cup \Sigma$}
$$
has a unique solution $(u_1, u_2)=\mathcal{T}_\Sigma(f_1,f_2)$, where $\mathcal{T}_\Sigma: \mathcal{C}_\Sigma(\overline\Omega)
\mapsto \mathcal{C}_\Sigma(\overline\Omega)
$
is a linear, compact and strongly positive operator, see~\cite{Suarezetal}. Moreover,
we denote by $\mathcal{T}_{\mathcal{N}}:=(-\Delta+\rho)^{-1}$, the inverse of $-\Delta+\rho$ in $\Omega$ under
homogeneous Neumann boundary conditions. Thus, let us define the inverse operator $\mathcal{T}$ as
$$
\mathcal{T}:=
(\mathcal{T}_\Sigma, \mathcal{T}_{\mathcal{N}})\mathbb{I}_3,
$$
with $\mathbb{I}_3$ the $3\times3$ identity matrix.

Consider the Banach space $\mathcal{U}:=\mathcal{C}_\Sigma(\overline\Omega) \times C(\overline\Omega)$
and let us denote ${\bf w}=(u_1,u_2,v)\in\mathcal{U}$.
Consider also the operator
$\mf{F}: \mathbb{R} \times \mathcal{U} \longrightarrow
\mathcal{U},
$
defined by
$$
  \mf{F}(\mu,{\bf w}):= {\bf w}- \mathcal{T}[(f_1,f_2,g)]\mathbb{I}_3{\bf w}
  $$
  where
$$
f_i=\lambda_i+\rho -\alpha_i u_{i}-a_i v ,\text{\ and\ } g=\mu+\rho- v-b_{1}u_{1}\chi_{\Omega_1}-b_{2}u_{2}\chi_{\Omega_2}.
$$
It is clear that ${\bf w}$ is a solution of \eqref{eq:main.system} if and only if $\mf{F}(\mu,{\bf w})=0$.

We know that $\mf{F}$ is of class
${C}^1$ and,
by elliptic regularity, $\mf{F}(\mu,\cdot,\cdot)$
is a compact perturbation of the identity for every $\mu\in\mathbb{R}$.
Moreover, $\mf{F}(\mu,{\bf S_\mu})=0$ for all $\mu\in\mathbb{R}$.
%
%
%
%
Differentiating  $\mf{F}$ with respect to ${\bf w}$, we have that,
$$
  D_{{\bf w}}  \mf{F}(\mu, {\bf w})\Big|_{{\bf w}={\bf S_\mu}}   =
  (\mathbb{I}_3-  \mathcal{T} M(x)),
$$
where
$$
M(x)=\left(\begin{array}{ccc}
   \lambda_1+\rho- 2\theta_{\l_1}\alpha_1 & 0& -a_1\theta_{\l_1} \\ 0 & \lambda_2+\rho - 2\theta_{\l_2}\alpha_2  & -a_2\theta_{\l_2}\ \\
   0 & 0 & \mu+\rho -b_{1}\theta_{\l_1}\chi_{\Omega_1}-b_{2}\theta_{\l_2}\chi_{\Omega_2}
   \end{array}\right).
$$
It is clear that, being a  compact perturbation of the
identity map, see~\cite{Bre}, the linear operator $D_{{\bf w}} \mf{F}(\mu,{\bf S_\mu})$  is a Fredholm
operator of index zero.
Moreover, we claim that it is injective, and hence  a linear topological isomorphism.
Indeed,  we denote
$$
\mf{D}_0:= D_{{\bf w}} \mf{F}(\mu_0,{\bf S_\mu}).
$$
First, observe that if $\overline{\bf w}=(\overline{u}_1,\overline{u}_2,\overline{v})\in {\rm Ker}(\mf{D}_0)$  then $\overline{v}$ is a eigenfunction associated with $\mu_0$, that is,
$$
-\Delta\overline{v}+ (b_{1}\theta_{\l_1}\chi_{\Omega_1}+b_{2}\theta_{\l_2}\chi_{\Omega_2})\overline{v}=\mu_0\overline{v}\quad\mbox{in $\O$,}\quad \partial_{\bf n}\overline{v}=0\quad\mbox{on $\partial\O$,}
$$
which provides us with the singular value $\mu_0$ given by \eqref{con:mu_zero}. Also, $\overline{u}_i$ verifies
$$
-\Delta \overline{u}_i+(2\theta_{\l_i}\alpha_i-\l_i)\overline{u}_i=-a_i\theta_{\l_i}\varphi_1\quad\mbox{in $\O_i$, }\quad \mathcal{I}(\overline{\bf u})=0\quad\mbox{on $\G\cup \Sigma$.}
$$
and thanks to (\ref{eigenposi}), there exists a unique solution of those two equations, see~\cite{Suarezetal}. Hence,
$$
{\rm Ker}(\mf{D}_0)={\rm Span}[{\bf \overline w}].
$$
On the other hand,
$$
\mf{D}_1 := D_\mu \big(D_{{\bf w}} \mf{F}(\mu,{\bf S_\mu})\big)\Big|_{\mu=\mu_0}= - \mathcal{T} \left(\begin{array}{ccc} 0 & 0 &0  \\ 0 &0 &0 \\ 0& 0& 1\end{array}\right).
$$
Consequently,
the following transversality condition holds
\begin{equation}
\label{trans_cond}
  \mf{D}_1 {\bf \overline w}\not\in {\rm Range}[  \mf{D}_0 ]
\end{equation} Indeed, arguing by contradiction we
suppose
that  there exists ${\bf w}\in \mathcal{U}$ such that
$$
 \mf{D}_0({\bf w})= \mf{D}_1 {\bf \overline w}.
$$
Now, considering only the third equation given by operator $\mf{D}_0$, we get that
$$
-\Delta v-(\mu_0 -b_{1}\theta_{\l_1}\chi_{\Omega_1}-b_{2}\theta_{\l_2}\chi_{\Omega_2})v=-\overline v.
$$
Observe that the other two equations provide us with equalities.
Then, multiplying such an third equation by $\overline v$ and integrating in $\Omega$, we arrive at
$$
\int_\Omega\overline v^2=0,
$$
which it is a contradiction.
Therefore, condition \eqref{trans_cond} holds, in other words,
$$ \mf{D}_1({\rm Ker}(\mf{D}_0))\oplus  {\rm Range}[  \mf{D}_0 ]= \mathcal{U}.$$
Consequently, according
to the main theorem of Crandall and Rabinowitz~\cite{CR}, see also \cite{julianlibro},  $(\mu_0, {\bf S_\mu})$ is a bifurcation point from
$(\mu, {\bf S_\mu})$ to a smooth curve of positive solutions of
\eqref{eq:main.system}.
\end{proof}

The following result shows the branch of coexistence states bifurcating from a value of the parameter $\mu$, once the other two main parameters $\lambda_1$ and $\lambda_2$ are fixed,
exists up to another value of the parameter $\mu$, providing us with the existence of positive solutions between those two values of $\mu$.  To this end, we denote the space
$$
\mathcal{E}:=\mathcal{C}^1_\Sigma(\overline\O)\times C^1(\overline\Omega).
$$
\begin{theorem}
\label{Theo_Bif2}
Let  $\lambda_i$ be fixed such that $\Lambda_1(-\l_1,-\l_2)<0$. Then, there exists a bounded continuum of solutions $\mathcal{C} \subset \mathbb{R}\times \mathcal{E}$ emanating from the semitrivial solution ${\bf S_\mu}$ at the value $\mu=\mu_0$ going up to the value of the parameter $\mu=\mu_1$, wich is given by
$$
\Lambda_1(-\lambda_1+a_1 {\mu_1},-\lambda_2 + a_2 {\mu_1})=0.
$$
As  a consequence,  there exists at least a  coexistence state of \eqref{eq:main.system}  for
\begin{equation}
\label{condi}
\mu\in (\min\{\mu_0,\mu_1\},\max\{\mu_0,\mu_1\}).
\end{equation}
\end{theorem}

\begin{proof} Thanks to Theorem~\ref{Theo_Bif} it is already proved that $\mu_0$ exists. Now, applying Theorem 7.2.2 of \cite{julianlibro}, and following the argument shown in Theorem 1.2 of \cite{cintra} (see also \cite{shi}) we conclude that a continuum $\mathcal{C} \subset \mathbb{R}\times \mathcal{E}$ of coexistence states emanates
from $\mu\equiv \mu_0$ such that, one of the following alternatives hold:
\begin{enumerate}
\item   $\mathcal{C}$ is unbounded in $ \mathbb{R}\times \mathcal{E}$;
\item there exists a sequence $(\mu_n,u_{1,n},u_{2,n},v_n)\in \mathcal{C}$  such that
$$(\mu_n,u_{1,n},u_{2,n},v_n)\to (\mu_\infty,0,0,0)\in \overline{ \mathcal{C}}
$$
for some $\mu_{\infty}\in \mathbb{R}$;
\item there exists a sequence $(\mu_n,u_{1,n},u_{2,n},v_n)\in \mathcal{C}$  such that
$$(\mu_n,u_{1,n},u_{2,n},v_n))\to (\mu_\infty,U_1,U_2,0)\in \overline{\mathcal{C}}
$$
for some $\mu_{\infty}\in \mathbb{R}$, $\mu_\infty\neq \mu_0$, and $U_i>0$ in $\O_i$.
\item there exists a sequence $(\mu_n,u_{1,n},u_{2,n},v_n)\in \mathcal{C}$  such that
$$
(\mu_n,u_{1,n},u_{2,n},v_n)\to (\mu_\infty,0,0,v_\infty)\in \overline{\mathcal{C}}
$$
for some $\mu_{\infty}\in \mathbb{R}$ and $v_\infty>0$ in $\Omega$.
\end{enumerate}
First, observe that thanks to Theorem \ref{cotassolu},
 we conclude
\begin{equation}
\label{eq:projection}
\mbox{Proj}_\mathbb{R}(\mathcal{C})\subset[0,\mu_*(\l_1,\l_2)],
\end{equation}
where $\mbox{Proj}_\mathbb{R}(\mu,u_1,u_2,v)=\mu$. Now, thanks to Theorem \ref{cotassolu}, we already have a priori bounds for the positive solution of~\eqref{eq:main.system}. Moreover,
by elliptic regularity, we find that they are also bounded in  $\mathcal{E}$. Hence,  $\mathcal{C}$ is bounded and we can conclude that alternative (1) is not possible.

Next, we assume that alternative (2) occurs. Thus, we have that
$$
0=\Lambda_1(-\l_1+a_1 v_n+\a_1u_{1,n},-\l_2+a_2v_n+\a_2 u_{2,n})\to \Lambda_1(-\l_1,-\l_2),
$$
arriving at a contradiction, by construction, since $\lambda_1$ and $\lambda_2$ are fixed with $\Lambda_1(-\l_1,-\l_2)<0$.

On the other hand, if alternative (3) occurs, then we have actually that $U_i=\theta_{\lambda_i}$ and, then, since by construction
$$
\mu_n=\sigma_1^\Omega[-\Delta + v_n+b_1 u_{1,n}\chi_{\Omega_1}+  b_2 u_{2,n}\chi_{\Omega_2};\mathcal{N}],
$$
passing to the limit, it yields to
$$
\mu_\infty=\sigma_1^\Omega[-\Delta + b_1 \theta_{\lambda_1}\chi_{\Omega_1}+  b_2 \theta_{\lambda_2} \chi_{\Omega_2};\mathcal{N}]=\mu_0,
$$
which it is again a contradiction due to $\mu_0\neq \mu_\infty$.
Consequently, it follows that alternative (4) occurs.
Indeed, observe that, by elliptic regularity $v_\infty$ is the positive solution of problem
$$
-\Delta v=v(\mu_\infty-v)\quad\mbox{in $\Omega$,}\quad  \partial_{\bf n} v=0\quad\mbox{on $\partial \Omega$,}
$$
and, hence, $v_\infty=\mu_\infty$. Moreover,
$$
0=\Lambda_1(-\l_1+a_1 v_n+\a_1u_{1,n},-\l_2+a_2v_n+\a_2 u_{2,n})\to\Lambda_1(-\l_1+a_1 \mu_\infty,-\l_2+a_2\mu_\infty)
$$
and then, by definition we find that  $\mu_\infty=\mu_1$.
Therefore, we actually have a continuum of positive solutions emanating from $\mu_0$ up to $\mu_1$.
\end{proof}


\section{Study of the region of coexistence}\label{sect:coexistence}

In this section we analyse the coexistence region given by Theorem \ref{Theo_Bif2}. To this aim, we will use the function $G(\lambda_1,\lambda_2)$ denoted by \eqref{eq:G_lambdas}.
Recall that this definition implies that $\mu_0=G(\l_1,\l_2)$, see~\eqref{con:mu_zero} and it is given by the definitions of the eigenvalue problems shown in Section 2.

The following result shows condition (\ref{condi}) in an equivalent form, which will be more appropriate in the sequel.

\begin{proposition}
\label{regiones}
If
\begin{equation}
\label{condi2}
(\mu-G(\lambda_1,\lambda_2))\cdot\Lambda_1(-\lambda_1+a_1\mu,-\lambda_2+a_2\mu)<0,
\end{equation}
then there exists at least a coexistence state of \eqref{eq:main.system}.
\end{proposition}
\begin{proof}
Assume without loss of generality that $\mu_0<\mu_1$ in~\eqref{condi}. In this case, Theorem \ref{Theo_Bif2} gives the existence of coexistence state for
$$
\mu\in (\mu_0,\mu_1).
$$
It is clear, see  \eqref{con:mu_zero} and \eqref{eq:G_lambdas}, that $\mu_0=G(\lambda_1,\lambda_2)$, hence $\mu>\mu_0$ is equivalent to $\mu-G(\lambda_1,\lambda_2)>0$. On the other hand, if $\mu<\mu_1$, it yields
$$
0=\Lambda_1(-\lambda_1+a_1 {\mu_1},-\lambda_2 + a_2 {\mu_1})>\Lambda_1(-\lambda_1+a_1 {\mu},-\lambda_2 + a_2 {\mu}).
$$
Hence, $\mu\in (\mu_0,\mu_1)$ is equivalent to \eqref{condi2}.
Note that one might treat similarly the case $\mu_0>\mu_1$.
\end{proof}
Accordingly, in the sequel, our aim will be studying the behaviour of the factors in~\eqref{condi2}.
To deal with $\Lambda_1$ we consider the following particular eigenvalue problem, assuming different values in each subdomain $\Omega_i$
\begin{equation}
\label{eign_prob}
\left\{\begin{array} {l@{\quad}l} -\Delta \varphi_i =\nu_i \varphi_i,&  \hbox{in}\quad \Omega_i,\\
\mathcal{I}({\bf \Phi})=0,& \hbox{on}\quad \Sigma\cup \Gamma. \end{array}\right.
\end{equation}
It is clear that $\nu_i$ is a principal eigenvalue of (\ref{eign_prob}) if and only if
\begin{equation}
\label{ecuimpor}
\Lambda_1(-\nu_1,-\nu_2)=0,
\end{equation}
and, hence, equation (\ref{ecuimpor}) defines a curve in $\mathbb{R}^2$ that we now describe.

\begin{figure}[ht!]
	\centering
\includegraphics[scale=0.25]{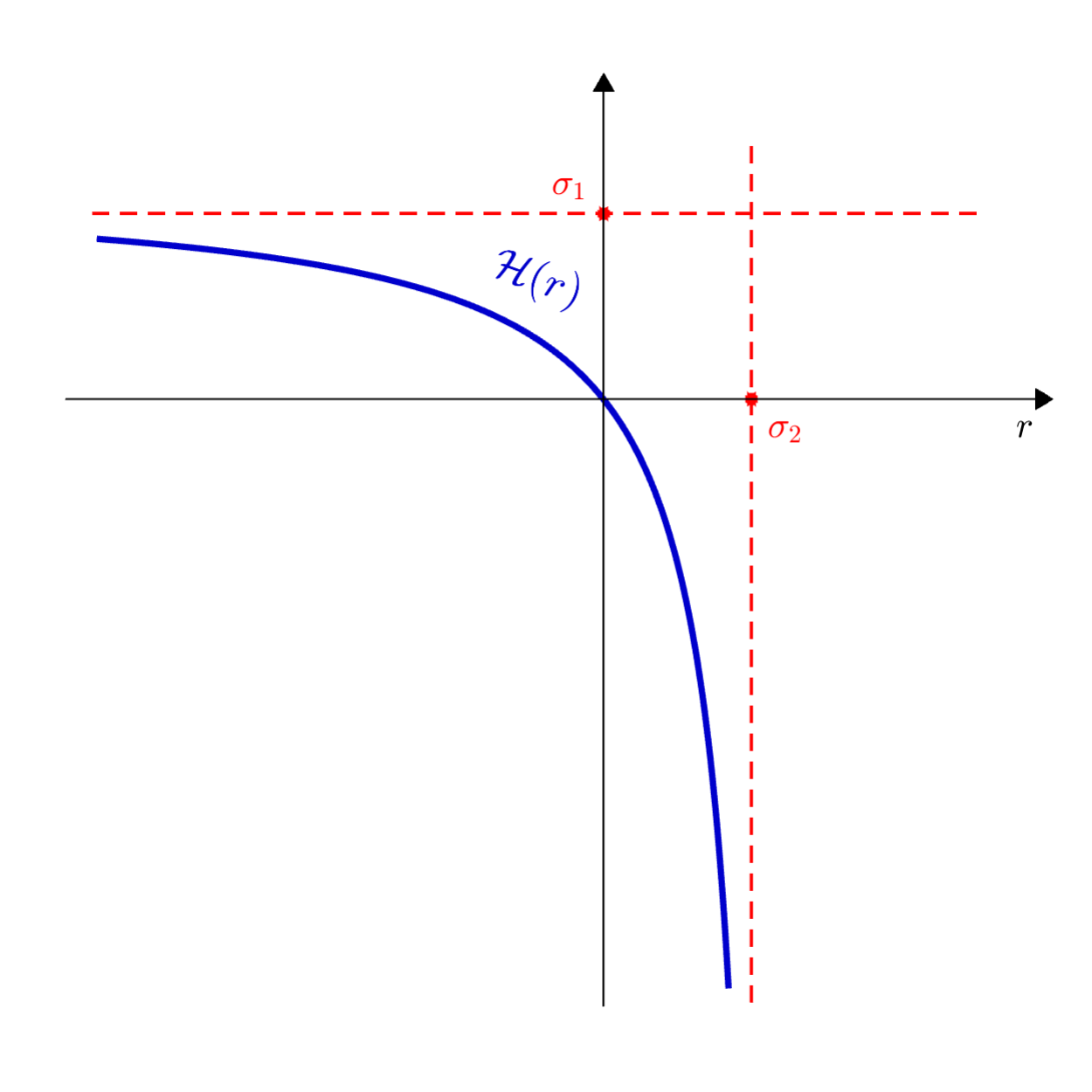}
\caption{Curve characterising~\eqref{ecuimpor}.}
\label{fig:func_h}
\end{figure}


\begin{lemma}
\label{autovalor}
There exists a continuous and decreasing map $
\mathcal{H}:(-\infty,\sigma_2)\mapsto \mathbb{R}
$, with $\mathcal{H}(0)=0$ and
$$
\lim_{r\to -\infty}\mathcal{H}(r)=\sigma_1,\qquad \lim_{r\to \sigma_2}\mathcal{H}(r)=-\infty,
$$
such that:
\begin{enumerate}
\item $\Lambda_1(-\nu_1,-\nu_2)=0$ if and only if $\nu_1=\mathcal{H}(\nu_2)$.
\item If $\nu_2\geq \sigma_2$, then $\Lambda_1(-\nu_1,-\nu_2)<0$ for all $\nu_1\in \mathbb{R}$.
\item If $\nu_2< \sigma_2$, then $\Lambda_1(-\nu_1,-\nu_2)<0 $ if $\nu_1> \mathcal{H}(\nu_2)$ (resp.$>0$, $\nu_1< \mathcal{H}(\nu_2)$).
\end{enumerate}
Moreover, $\mathcal{H}$ is differentiable and
\begin{equation}
\label{nova}
\mathcal{H}'(0)=\lim_{s\to 0}\frac{\mathcal{H}(s)}{s}=-\frac{\gamma_1}{\g_2}\frac{|\O_2|}{|\O_1|}.
\end{equation}
\end{lemma}
\begin{proof}
First note that with the exception of (\ref{nova}), this lemma follows from \cite{Suarezetal2}.

To prove (\ref{nova}) we observe that, by  (1) above
$$
\Lambda_1(-\D-\mathcal{H}(s),-\D -s)=0,
$$
and so, dividing by $s$, we get
\begin{equation}
\label{ceroese}
0=\Lambda_1\left(-\frac{1}{s}\D-\frac{\mathcal{H}(s)}{s},-\frac{1}{s}\D -1\right).
\end{equation}
Assume that $\frac{\mathcal{H}(s)}{s}\to L\in [-\infty,+\infty]$ as $s\to 0$.
 Indeed, integrating both equations of the eigenvalue problem associated with the principal eigenvalue \eqref{ceroese}, in the corresponding $\Omega_i$, after multiplying the second by $\gamma_1/\gamma_2$ we find that
 $$\frac{1}{s} \int_{\Omega_1} \Delta \varphi_1+  \int_{\Omega_1} \frac{\mathcal{H}(s)}{s} \varphi_1+ \frac{1}{s} \frac{\gamma_1}{\gamma_2} \int_{\Omega_2} \Delta \varphi_2 + \frac{\gamma_1}{\gamma_2} \int_{\Omega_2} \varphi_2=0.$$
 Now, due to Gauss Divergence Theorem and the boundary conditions  $\mathcal{I}({\bf \Phi})$ denoted by \eqref{eq:BC_KK} it follows that
 $$\int_{\Omega_1} \Delta \varphi_1 +  \frac{\gamma_1}{\gamma_2} \int_{\Omega_2} \Delta \varphi_2=0.$$
 Consequently,
 $$ \int_{\Omega_1} \frac{\mathcal{H}(s)}{s} \varphi_1=- \frac{\gamma_1}{\gamma_2} \int_{\Omega_2} \varphi_2,$$
 and since the eigenfunctions are bounded in  $H^1(\Omega_1)\times H^1(\Omega_2)$, then they converge weakly in $H^1(\Omega_1)\times H^1(\Omega_2)$ and strongly in $L^2(\Omega_1)\times L^2(\Omega_2)$,
 so that the assumption on the convergence of the term $\frac{\mathcal{H}(s)}{s}$ is now justified.

First, observe that from (\ref{ceroese}) we can deduce that $L\in (-\infty,+\infty)$. Moreover, by~\eqref{dainfinito}, if we pass to the limit in (\ref{ceroese}) as $s$ goes to zero, we arrive at
$$
0=\lim_{s\to 0} \Lambda_1\left(-\frac{1}{s}\Delta-\frac{\mathcal{H}(s)}{s},-\frac{1}{s}\Delta -1\right)=\frac{-\gamma_2|\Omega_1|L-\gamma_1|\Omega_2|}{\gamma_2|\Omega_1|+\gamma_1|\Omega_2|},
$$
and~\eqref{nova} is proved.
\end{proof}


Now, as a first step we study the sign of $\Lambda_1(-\lambda_1+a_1\mu,-\lambda_2+a_2\mu)$. Observe that, according to Lemma~\ref{autovalor}, and for the particular eigenvalue
we are analysing now, there will be a critical value, $\sigma_2= -\lambda_2+a_2\mu$,
that it will distinguish different behaviours of $\Lambda_1$. Hence, we define the map
\begin{equation}
  \label{functionG}
\mathcal{G}:\left(\frac{\l_2-\sigma_2}{a_2},+\infty\right)\mapsto \mathbb{R},\quad \mathcal{G}(\mu):=a_1\mu+\mathcal{H}(\l_2-a_2\mu).
\end{equation}
It is straightforward to see, using Lemma~\ref{autovalor}, that $\mathcal{G}$ is a continuous and increasing function. Moreover,
$$
\lim_{\mu\to \frac{\l_2-\sigma_2}{a_2}}\mathcal{G}(\mu)=-\infty,\quad \lim_{\mu\to+\infty}\mathcal{G}(\mu)=+\infty.
$$
%

\begin{lemma}
\label{rama1} Let $\mathcal{G}$ be given by~\eqref{functionG}.
\begin{enumerate}
	\item If  $\mu\leq (\l_2-\sigma_2)/a_2$, then,
$$
\Lambda_1(-\l_1+a_1\mu,-\l_2+a_2\mu)<0\quad\mbox{for all $\l_1\in\mathbb{R}$}.
$$
\item If  $\mu>(\l_2-\sigma_2)/a_2$, then,
$$
\Lambda_1(-\l_1+a_1\mu,-\l_2+a_2\mu)<0\quad\mbox{if and only if}\quad \l_1>\mathcal{G}(\mu),
$$
and
$$
\Lambda_1(-\l_1+a_1\mu,-\l_2+a_2\mu)>0\quad\mbox{if and only if}\quad \l_1<\mathcal{G}(\mu).
$$
\end{enumerate}
\end{lemma}
\begin{proof}
Note that $
\Lambda_1(-\l_1+a_1\mu,-\l_2+a_2\mu)=\Lambda_1(-(\l_1-a_1\mu),-(\l_2-a_2\mu))
$. Then, according to Lemma~\ref{autovalor} we obtain the following:

If $\l_2-a_2\mu\geq \sigma_2$, that is $\mu\leq (\l_2-\sigma_2)/a_2$, then 
$$
\Lambda_1(-\l_1+a_1\mu,-\l_2+a_2\mu)<0\quad\mbox{for all $\l_1\in\mathbb{R}$}.
$$
On the contrary, if $\l_2-a_2\mu<\sigma_2$, then, 
$$\Lambda_1(-\l_1+a_1\mu,-\l_2+a_2\mu)<0\quad\mbox{if and only if
$\l_1-a_1\mu>\mathcal{H}(\l_2-a_2\mu)$,}
$$
or equivalently, $\l_1>\mathcal{G}(\mu)$.
Analogously, $\Lambda_1(-\l_1+a_1\mu,-\l_2+a_2\mu)>0$ if and only if $\l_1<\mathcal{G}(\mu).$
\end{proof}
Next, we focus on the analysis of $G(\lambda_1,\lambda_2)$ defined in \eqref{eq:G_lambdas}.
To this aim, we shall assume that one of the parameters moves, let us say $\lambda_1$,  while the other one, $\lambda_2$, is fixed. Thus, according to the
eigenvalue problems stated by \eqref{omega2} and the definition of the eigenvalues \eqref{eq:def.sigmas}, fixing the parameter $\lambda_2$ we define
\begin{equation}
\label{def:eig_mu}
g(\lambda_1):=\sigma_1^\Omega[-\Delta + b_1 \theta_{\lambda_1}\chi_{\Omega_1}+  b_2 \theta_{\lambda_2} \chi_{\Omega_2};\mathcal{N}].
\end{equation}
Recall also the definition of the large solution $L_{\lambda_2}$, see~\eqref{eq:large}.

\begin{lemma}
\label{lemma.glambda1}
Let $\lambda_2$ be fixed. Then, the map $\l_1\mapsto g(\l_1)$ is continuous, increasing and
$$
\lim_{\lambda_1 \to +\infty} g(\lambda_1)= \sigma_1^{\Omega_2}[-\Delta + b_2 L_{\lambda_2};\mathcal{D};\mathcal{N}]. 
$$
where such a limiting eigenvalue is the principal eigenvalue of problem
\begin{equation}
\label{eq:classic_Dir_Neu}
 \left\{\begin{array} {l@{\quad}l} (-\Delta + b_2 L_{\lambda_2}) \varphi =\sigma \varphi,& \hbox{in}\quad \Omega_2,\\
 \varphi=0,& \hbox{on}\quad \Sigma,\\
 \partial_{\bf n} \varphi=0, & \hbox{on}\quad \Gamma.
\end{array}\right.
\end{equation}
with homogeneous Dirichlet boundary conditions on $\Sigma$ and homogenous Neumann boundary conditions on $\Gamma$, such that $\partial \Omega_2=\Sigma \cup \Gamma$.
Moreover:
\begin{enumerate}
\item  If $\lambda_2\geq \sigma_2$, then $g$ is defined for all $\l_1\in \mathbb{R}$ and
\begin{equation}
\label{limiminus}
\lim_{\l_1\to-\infty}g(\lambda_1)= \sigma_1^{\O_2}[-\Delta + b_2 \omega_{2};\mathcal{N}].
\end{equation}
\item If $\lambda_2< \sigma_2$, then $g$ is defined for  $\l_1\in [\mathcal{H}(\l_2),\infty)$ and $g(\mathcal{H}(\l_2))=0$.
\end{enumerate}
\end{lemma}
\begin{proof}
It is clear that $\l_1\mapsto g(\l_1)$ is continuous, increasing and, also, thanks to the monotonicity properties of the first eigenvalues (see~\cite{Santi} for further details)
$$
g(\l_1)<\sigma_1^{\Omega_2}[-\Delta + b_1 \theta_{\lambda_1}\chi_{\Omega_1}+  b_2 \theta_{\lambda_2} \chi_{\Omega_2};\mathcal{D};\mathcal{N}]=\sigma_1^{\Omega_2}[-\Delta+ b_2 \theta_{\lambda_2} \chi_{\Omega_2};\mathcal{D};\mathcal{N}].
$$
Then, $g(\lambda_1)$ is bounded above. Thus, there exists $g^*<\infty$ such that
$$
\lim_{\l_1\to +\infty}g(\l_1)=g^*.
$$
Now, to identify such an specific value $g^*$, take a sequence $\l_{1,n}\to +\infty$ and consider $\varphi_{1,n}$ the  positive eigenfunction associated with $g(\l_{1,n})$ such that $\|\varphi_{1,n}\|_2=1$. Hence, we can conclude that
$$
\int_\Omega |\nabla \varphi_{1,n}|^2\leq C,
$$
and, then, applying the compact embedding of $H^1(\Omega)$ into $L^2(\Omega)$ we can extract a subsequence, again labelled  $\{\varphi_{1,n}\}$, converging, weakly in $H^1(\Omega)$ and strongly in $L^2(\Omega)$, i.e.
$$
 \varphi_{1,n}\rightharpoonup\varphi_{1,\infty}\geq 0 \quad\mbox{in $H^1(\Omega)$},\quad  \varphi_{1,n}\to\varphi_{1,\infty} \quad\mbox{in $L^2(\Omega)$ with $\|\varphi_{1,\infty}\|_2=1$.}
$$
Furthermore, we claim that
\begin{equation}
\label{ho1}
\varphi_{1,\infty}\in H_{0,\Sigma}^1(\Omega_2):=\{u\in H^1(\Omega_2):\mbox{$u=0$ on $\Sigma$}\}.
\end{equation}
To prove it we argue by contradiction. Indeed, assume that $\varphi_{1,\infty}>0$ in $D$, for some $D\subset \Omega\setminus \Omega_2$ and take $\psi>0$ in $D$
such that $\psi\in C_0^\infty(D)$. Then, using the definition of $\varphi_{1,n}$,  as  the  positive eigenfunction associated with $g(\l_{1,n})$, multiplying the equation,
of its corresponding eigenvalue problem, by $\psi$ and integrating by parts, we conclude that
$$
-\int_D \varphi_{1,n} \Delta \psi+b_1\int_D \theta_{\lambda_{1,n}} \varphi_{1,n} \psi=g(\lambda_{1,n})\int_D  \varphi_{1,n}\psi.
$$
Since $\theta_{\lambda_{1,n}}\to\infty$, when $\lambda_{1,n}\to +\infty$
(see Proposition\;\ref{Prop_v0} proved in~\cite{AlvarezBrandleMolinaSuarez} for further details)
we arrive at a contradiction, because
$$
\int_D \theta_{\lambda_{1,n}} \varphi_{1,n}\psi\to \infty,\quad\text{and}\quad \int_D \varphi_{1,\infty} \Delta \psi+g^*\int_D  \varphi_{1,\infty}\psi\leq C.
$$
Hence, we conclude that $\varphi_{1,\infty}=0$ in $D$, which implies (\ref{ho1}).

Subsequently, take now the test function, $\phi$ so that 
$\phi\in C_{c}^\infty(\O_2\cup \Gamma)$ and  $\phi=0$ in $\Omega_1$. Then, according to the definition of $g(\lambda_1)$ \eqref{def:eig_mu} in terms of
the eigenvalue problem \eqref{omega2}, multiplying such a corresponding equation by $\phi$ and integrating again by parts, we have that
$$
\int_{\Omega_2}\nabla \varphi_{1,n}\cdot\nabla \phi+b_2\int_{\Omega_2}\theta_{\lambda_2}\varphi_{1,n} \phi=g(\lambda_{1,n})\int_{\Omega_2}\varphi_{1,n} \phi.
$$
Passing to the limit as $\lambda_{1,n}\to\infty$ it yields to
$$
\int_{\Omega_2}\nabla \varphi_{1,\infty}\cdot\nabla \phi+b_2\int_{\Omega_2}L_{\lambda_2}\varphi_{1,\infty} \phi=g^*\int_{\Omega_2}\varphi_{1,\infty} \phi,\qquad\forall \phi\in C_{c}^\infty(\O_2\cup \Gamma),
$$
where we have used that $\theta_{\lambda_2}$ tends to the large solution $L_{\lambda_2}$ given by problem \eqref{eq:large} as $\lambda_{1,n}$ goes to infinity;
see again~\cite{AlvarezBrandleMolinaSuarez} for further details where Proposition\;\ref{Prop_v0} is proved. Consequently,  since $\varphi_{1,\infty}\geq 0$, $\varphi_{1,\infty}\neq 0$, we can finally conclude that
$$
g^*=\sigma_1^{\Omega_2}[-\Delta+ b_2 L_{\lambda_2};\mathcal{D};\mathcal{N}],
$$
as the principal eigenvalue of problem \eqref{eq:classic_Dir_Neu} (see Lemma 3.1 in \cite{du} and Lemma 8.3 in \cite{julianhouston} for eigenvalue problems similar to \eqref{eq:classic_Dir_Neu}).

Next, we prove the final part of the lemma:

First assume that $\l_2\geq\sigma_2$. Thanks to the results obtained in~\cite[Section 5]{AlvarezBrandleMolinaSuarez}, for such a particular case we know that $(\theta_{\lambda_1},\theta_{\lambda_2})$ exists
for all $\lambda_1\in \mathbb{R}$ and due to \cite[Proposition 5.7]{AlvarezBrandleMolinaSuarez} as $\lambda_1\to-\infty$
$$
\theta_{\l_1}\to 0\quad \hbox{and}\quad
\theta_{\l_2}\to \omega_2.
$$
Therefore, applying again the limiting argument we performed above in terms of $\sigma_1^{\Omega}$, see~\cite{AlvarezBrandleMolinaSuarez} again, we get (\ref{limiminus}).

 On the other hand, if $\lambda_2<\sigma_2$ we get from Lemma~\ref{autovalor} that $(\theta_{\l_1},\theta_{\l_2})$ exists if $\l_1>\mathcal{H}(\l_2)$ and by definition $g(\mathcal{H}(\l_2))=0$.
 \end{proof}
\begin{figure}[ht!]
	\centering
\includegraphics[scale=0.3]{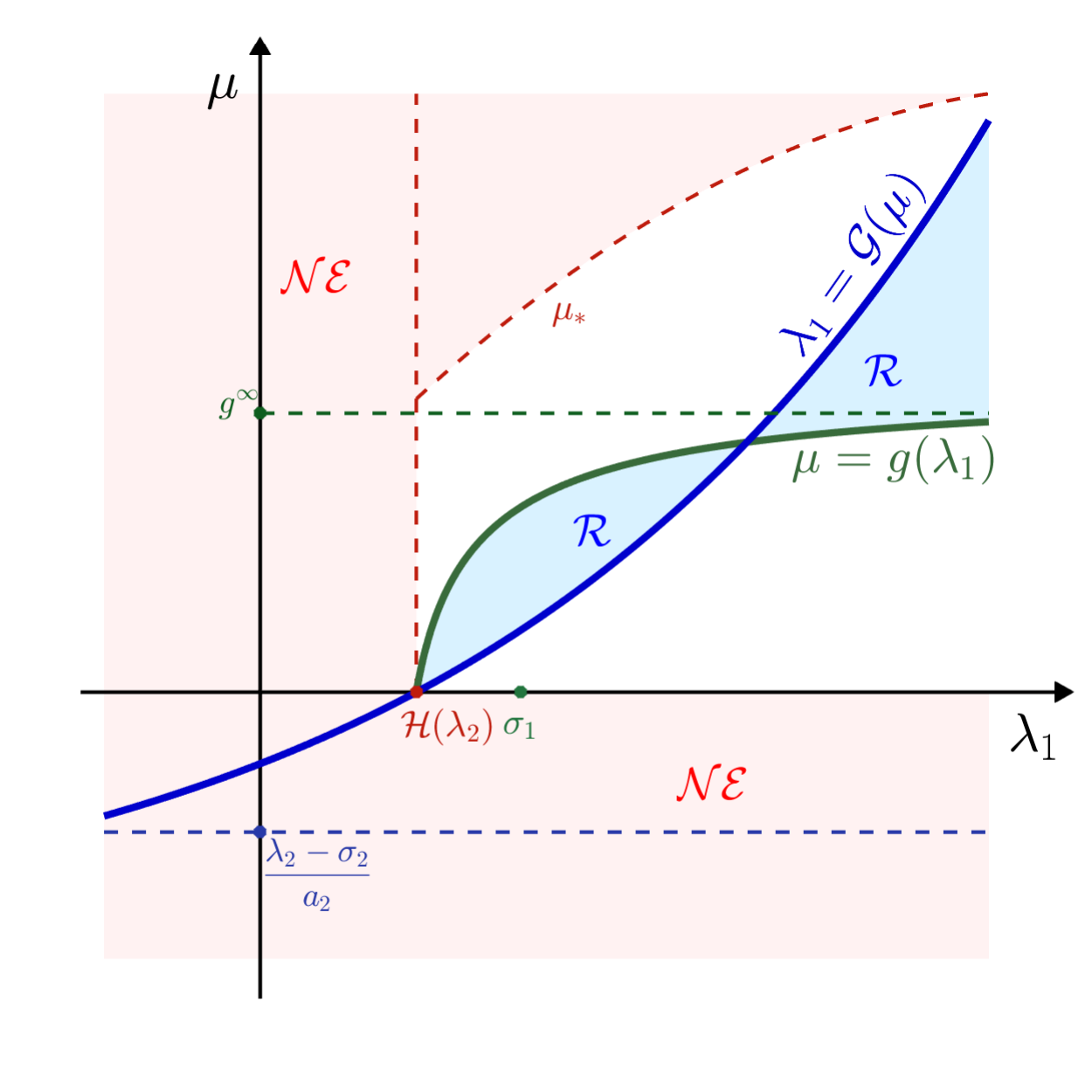}
\includegraphics[scale=0.3]{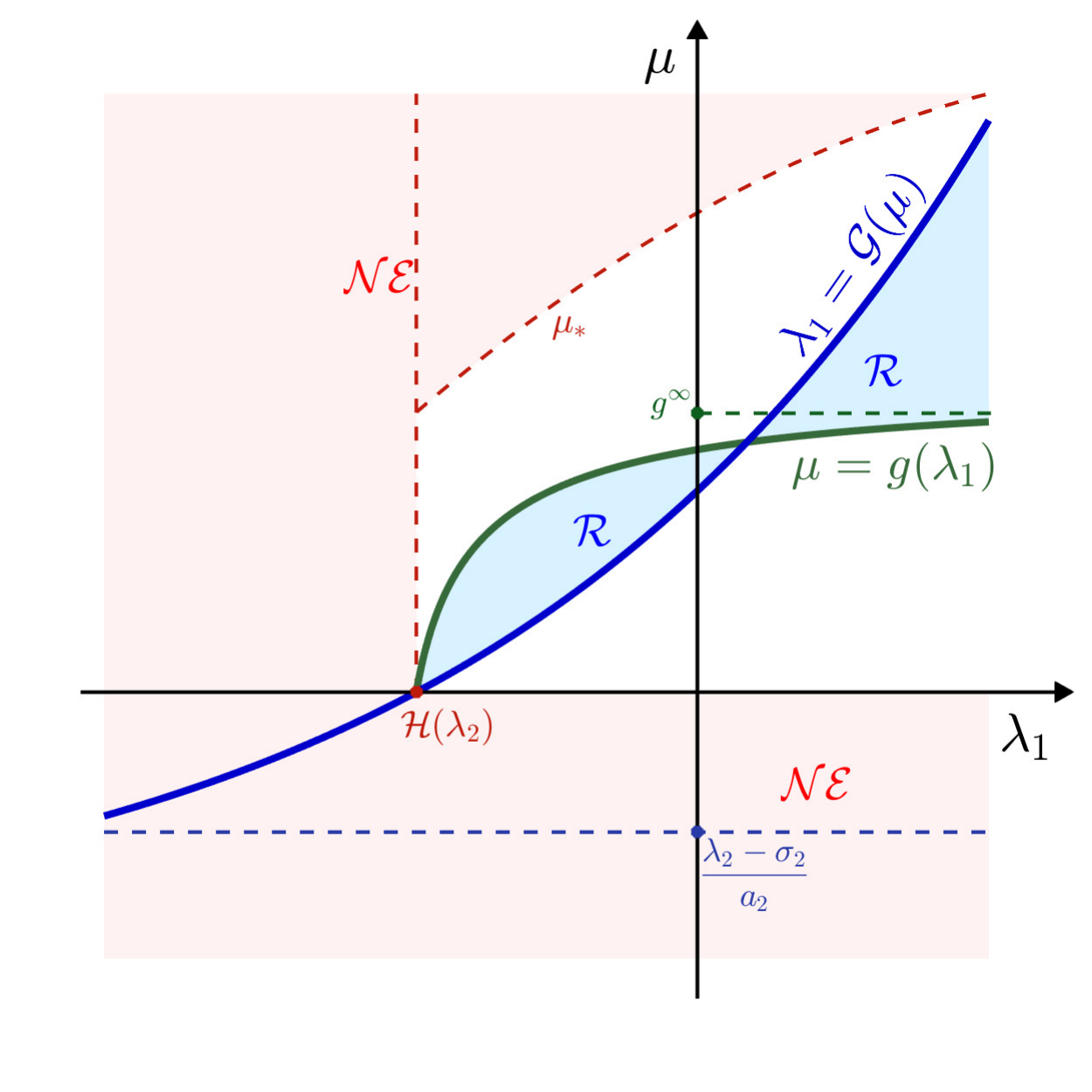}
\caption{Case $\l_2<\sigma_2$: In this case, $\Lambda_1(-\l_1,-\l_2)<0$ is equivalent to $\l_1>\mathcal{H}(\l_2).$ Hence, $g(\l_1)$ is defined for $\l_1\geq \mathcal{H}(\l_2)$ and $g^\infty= \lim_{\l_1\to +\infty}g(\l_1)$. Observe that $\mathcal{G}(0)=\mathcal{H}(\l_2)$. Left: case $\mathcal{H}(\l_2)>0$, that is $\l_2<0$; Right: case $\mathcal{H}(\l_2)<0$, that is $\l_2>0$. Here $\mathcal{R}$ denotes the coexistence states region and $\mathcal{NE}$ denotes the non-existence region defined by Corollary \ref{col:nece_cond} and Proposition \ref{nomularge}.}
\label{fig1}
\end{figure}
\begin{figure}[ht!]
	\centering
\includegraphics[scale=0.27]{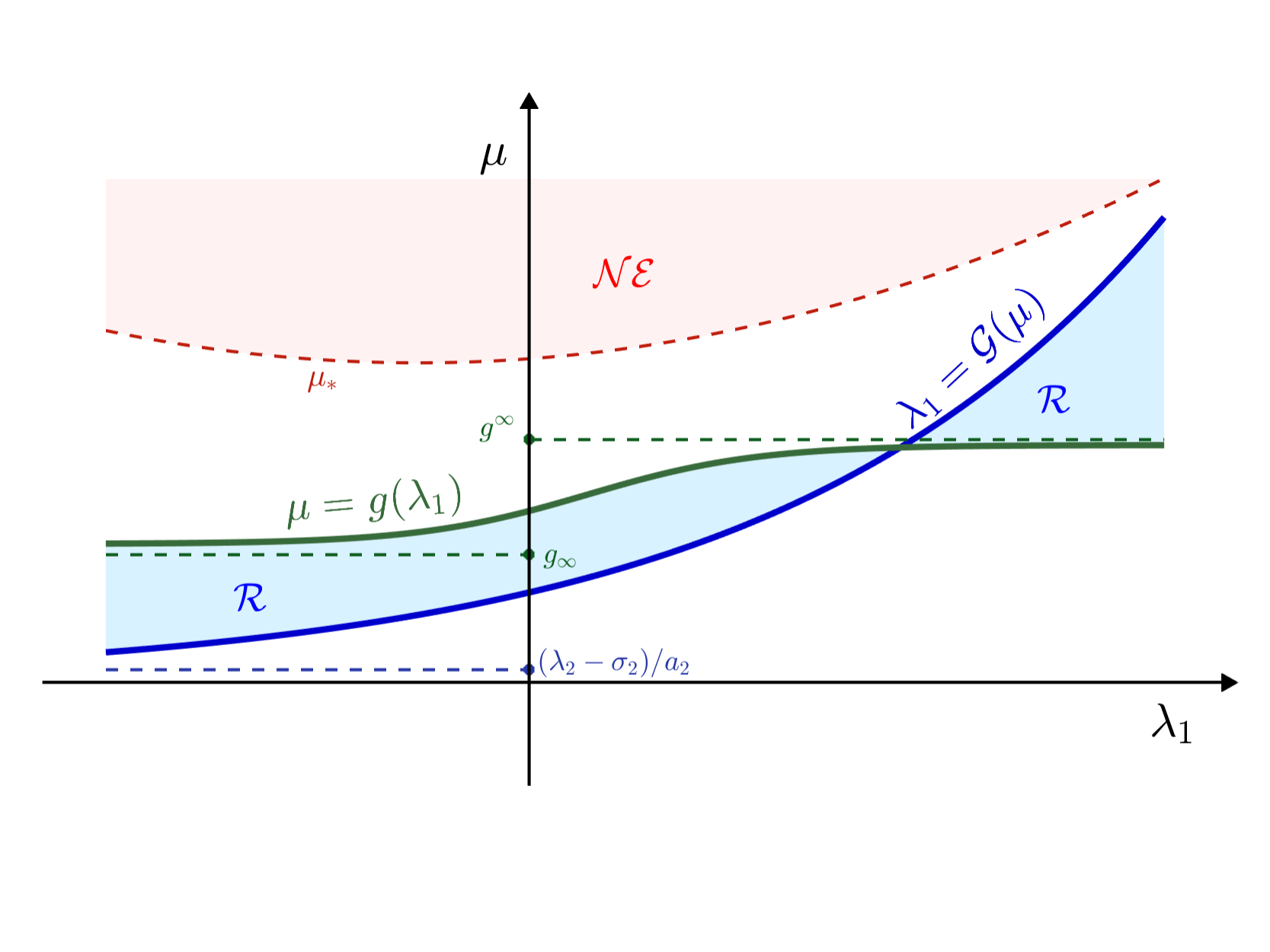}
\includegraphics[scale=0.27]{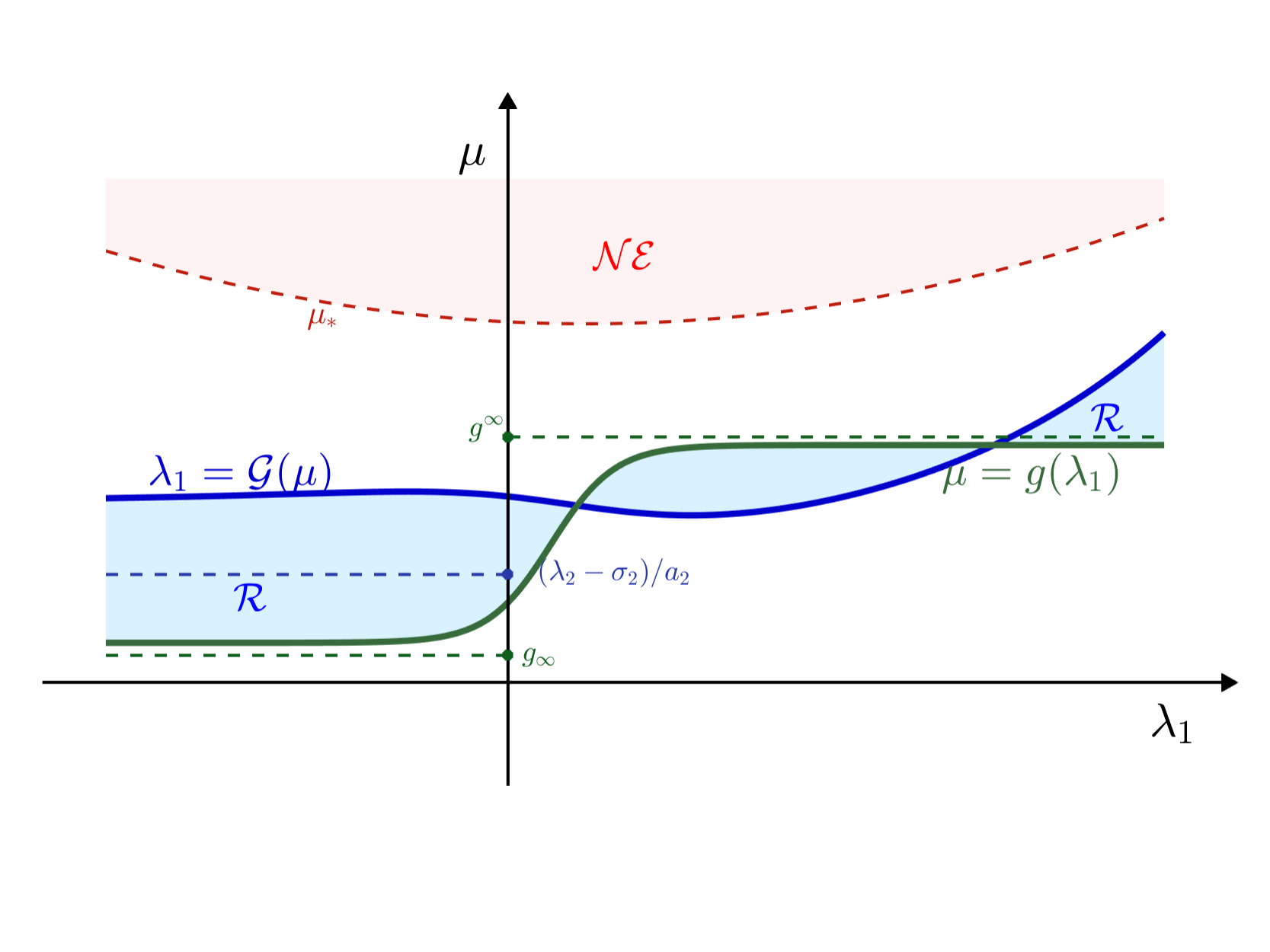}
\caption{Case $\l_2\geq \sigma_2$: In this case, $\Lambda_1(-\l_1,-\l_2)<0$ for all $\l_1\in\mathbb{R}$, and then $g(\l_1)$ is defined for all $\l_1\in\mathbb{R}$. Here, $g_\infty=\lim_{\l_1\to -\infty}g(\l_1)$ and $g^\infty= \lim_{\l_1\to +\infty}g(\l_1)$. Left: case $(\l_2-\sigma_2)/a_2<g_\infty$; Right: case $(\l_2-\sigma_2)/a_2>g_\infty$. Again $\mathcal{R}$ and $\mathcal{NE}$ denote the coexistence states and the non-existence  regions, respectively.}
\label{fig2}
\end{figure}

Throughout the rest of this section we will assume that  $\l=\l_1=\l_2$. Thus, under such an assumption it is possible to obtain sharper results than the previous Lemmas~\ref{rama1} and~\ref{lemma.glambda1}. Let us
 denote by $ (\theta^1_{\lambda}, \theta^2_{\lambda})$ the solution of~\eqref{eq:membrane.logis} with $\l=\l_1=\l_2$ and rewrite condition~\eqref{condi2}, for the particular problem in hand, as
\begin{equation}
\label{condi2igual}
(\mu-g(\l))\cdot\Lambda_1(-\l+a_1\mu,-\l+a_2\mu)<0, 
\end{equation}
with
$$
g(\l):=G(\l,\l)=\sigma_1^\Omega[-\Delta + b_1 \theta^1_{\lambda}\chi_{\Omega_1}+  b_2 \theta^2_{\lambda} \chi_{\Omega_2};\mathcal{N}].
$$

Recall that  $ (\theta^1_{\lambda}, \theta^2_{\lambda})$ is positive if and only if $\l>0$, just applying a comparison argument as the one used to provide a necessary condition for the existence of solutions shown in Corollary\;\ref{col:nece_cond}.
\begin{lemma}
\label{igualu}
The map $g:[0,\infty)\mapsto [0,\infty)$ is continuous, increasing and
\begin{equation}
\label{limi1}
\lim_{\l\to 0^+}g(\l)=0,\qquad \lim_{\l\to+\infty} g(\l)=+\infty.
\end{equation}
Moreover,
\begin{equation}
\label{limi2}
\lim_{\l\to +\infty}\frac{g(\l)}{\l}=\min\left\{\frac{b_1}{\alpha_1},\frac{b_2}{\alpha_2}\right\}\text{\; and \;} \lim_{\l\to 0^+}\frac{g(\l)}{\l}=\frac{(\g_1|\O_2|+\g_2|\O_1|)(b_1|\O_1|+b_2|\O_2|)}{\g_1\a_2|\O_2|+\g_2\a_1|\O_1|}.
\end{equation}
\end{lemma}
\begin{proof}
According to Proposition\;\ref{exislogismem} the map $\l\mapsto (\theta^1_{\lambda}, \theta^2_{\lambda})$ is increasing and continuous hence, we conclude, by monotonicity of the first eigenvalue with respect to the potential, that $g(\l)$ is also continuous and increasing.

Moreover, due to (\ref{cotaimpor}) it follows that $\|\theta^i_{\lambda}\|_\infty\to 0$   as $\l\to 0^+$ so that one can obtain that
$$
\lim_{\l\to 0^+}g(\l)=0.
$$
On the other hand,
applying again the fact that the map $\l\mapsto (\theta^1_{\lambda}, \theta^2_{\lambda})$ is increasing and continuous, and that
$$
\theta^i_{\lambda}\geq \omega_i\geq\frac{(\lambda-\sigma_i)}{\alpha_i}\varphi_1^i,\quad i=1,2,
$$
where $\omega_i$ are the solutions of (\ref{eq:logisuna}) and (\ref{omega2}) with $\l_i=\lambda$, respectively,  and $\varphi_1^i$ is the positive eigenfunction associated to $\sigma_i$ with $\|\varphi_1^i\|_\infty=1$  in $\Omega_i$ for $i=1,2$.  By definition of the function $g(\lambda)$ we have that
$$
g(\l)\geq \sigma_1^\Omega[-\Delta + b_1 \frac{(\lambda-\sigma_1)}{\alpha_1}\varphi_1^1\chi_{\Omega_1}+  b_2\frac{(\lambda-\sigma_2)}{\alpha_2}\varphi_1^2\chi_{\Omega_2};\mathcal{N}],
$$
and passing to the limit as $\l\to+\infty$ we prove (\ref{limi1}).

To prove~\eqref{limi2}, observe that
$$
\frac{g(\l)}{\l}=\sigma_1^\O[-\frac{1}{\l}\Delta + c_\l;\mathcal{N}], \quad \hbox{with}\quad c_\l:=b_1 \frac{\theta^1_{\lambda}}{\l}\chi_{\Omega_1}+  b_2 \frac{\theta^2_{\lambda}}{\l}\chi_{\Omega_2}.
$$
Since
$$
\lim_{\lambda\to+\infty} \frac{\theta_\l^i}{\l}\to \frac{1}{\alpha_i}, \quad\mbox{\quad and\quad\quad} \lim_{\l\to 0^+}\frac{\theta_\l^i}{\l}=\frac{\g_2|\O_1|+\g_1|\O_2|}{\g_2\alpha_1|\O_1|+\g_1\alpha_2|\O_2|},
$$
see~\cite[Corollary 5.5]{AlvarezBrandleMolinaSuarez}, we get, on one hand, that
\begin{equation}
\label{eq:clambda.infty}
\lim_{\l\to+\infty}c_\l= b_1\frac{1}{\alpha_1}\chi_{\Omega_1}+b_2 \frac{1}{\alpha_2}
\chi_{\Omega_2}\quad\mbox{in $L^\infty(\O)$,}
\end{equation}
and on the other
\begin{equation}
\label{eq:clambda.zero}
\lim_{\l\to 0^+} c_\l= \frac{\g_2|\O_1|+\g_1|\O_2|}{\g_2\alpha_1|\O_1|+\g_1\alpha_2|\O_2|}(b_1 \chi_{\Omega_1}+b_2 \chi_{\Omega_2})\quad\mbox{in $L^\infty(\O)$.}
\end{equation}
Finally, following~\cite[Proposition 2.2]{AlvarezBrandleMolinaSuarez} it can be shown that
$$
\lim_{\lambda\to +\infty}\sigma_1^\Omega[-\frac{1}{\lambda}\Delta+c_\lambda;{\mathcal{N}}]=\left(\lim_{\lambda\to+\infty} c_\lambda\right)_L=
\min\left\{\frac{b_1}{\alpha_1},\frac{b_2}{\alpha_2}\right\}.
$$
Also,
since $\sigma_1^D[-\Delta;{\mathcal{N}}]=0$, it follows by  \cite[Proposition 2.3]{AlvarezBrandleMolinaSuarez} that
$$
\lim_{\lambda\to 0^+ }\sigma_1^\Omega[-\frac{1}{\lambda}\Delta+c_\lambda;{\mathcal{N}}]=
\int_\Omega(\lim_{\lambda\to 0^+} c_\lambda)=  \frac{\g_2|\O_1|+\g_1|\O_2|}{\g_2\alpha_1|\O_1|+\g_1\alpha_2|\O_2|}(b_1 |\Omega_1|+b_2 |\Omega_2|).
$$
This concludes the proof.
\end{proof}

%

\begin{lemma}
\label{igual2}
For all $\mu \geq 0$, there exists a continuous and increasing function $\widehat{\mathcal{G}}:[0,+\infty)\mapsto [0,+\infty)$, such that
\begin{equation}
\label{senegal}
\widehat{\mathcal{G}}(\mu)\leq\sigma_2+ \min\{a_1,a_2\}\mu,
\end{equation}
 Moreover,
$$
\lim_{\mu\to 0^+}\widehat{\mathcal{G}}(\mu)=0,\qquad \lim_{\mu\to \infty}\widehat{\mathcal{G}}(\mu)=\infty,
$$
and
\begin{equation}
\label{oro}
\lim_{\mu\to 0^+}\frac{\widehat{\mathcal{G}}(\mu)}{\mu}=\frac{\g_1a_2|\O_2|+\g_2a_1|\O_1|}{\g_1|\O_2|+\g_2|\O_1|},\qquad\lim_{\mu\to \infty}\frac{\widehat{\mathcal{G}}(\mu)}{\mu}=\min\{a_1,a_2\}.
\end{equation}
Furthermore, 
$$
\Lambda_1(-\lambda+a_1\mu,-\lambda+a_2\mu)<0\quad\mbox{if and only if}\quad \lambda>\widehat{\mathcal{G}}(\mu),
$$
and
$$
\Lambda_1(-\lambda+a_1\mu,-\lambda+a_2\mu)>0\quad\mbox{if and only if}\quad \lambda<\widehat{\mathcal{G}}(\mu).
$$
\end{lemma}

\begin{proof}
Let us consider an auxiliary function ${h}:(-\infty,\sigma_2)\mapsto (-\infty,\infty)$ given by 
$$h(\sigma):=-\sigma +\mathcal{H}(\sigma)
$$ 
where the map $\mathcal{H}$ is defined in Lemma \ref{autovalor}. From Lemma~\ref{autovalor}, we have that $h$ is derivable and  decreasing, $h(0)=0$ and
$$
\lim_{\sigma\to -\infty}h(\sigma)=+\infty,\qquad\lim_{\sigma\to \sigma_2^-}h(\sigma)=-\infty.
$$
Then, for any $\mu\geq 0$ there exists a unique value $\sigma_0:=\sigma_0(\mu)$ such that
$$
h(\sigma_0(\mu))=(a_2-a_1)\mu.
$$
Observe that, from the Implicit Function Theorem, the map $\mu\mapsto \sigma_0(\mu)$ is derivable, $\sigma_0(\mu)<\sigma_2$ and the monotonicity of $\sigma_0(\mu)$ depends on whether $a_1>a_2$ or not,
$$
h'(\sigma_0(\mu))\sigma_0'(\mu)=(a_2-a_1).
$$
If $a_2>a_1$, then $\sigma_0(\mu)<0$, decreasing and
$$
\lim_{\mu \to 0}\sigma_0(\mu))=0, \quad \lim_{\mu \to +\infty}\sigma_0(\mu)=-\infty.
$$
On the other hand, if $a_2<a_1$, then $\sigma_0(\mu)>0$, increasing and
$$
\lim_{\mu \to 0}\sigma_0(\mu)=0, \quad \lim_{\mu \to +\infty}\sigma_0(\mu)=\s_2.
$$
Finally, $\sigma_0(\mu)=0$ for all $\mu$, if $a_1=a_2$.

Consequently, given $\mu\geq 0$ define
\begin{equation}
\label{defi1}
\widehat{\mathcal{G}}(\mu):=\sigma_0(\mu)+a_2\mu.
\end{equation}
Since, by definition of $h$,
\begin{equation}
\label{defi0}
 -\sigma_0(\mu)+ \mathcal{H}(\sigma_0(\mu))=(a_2-a_1)\mu,
\end{equation}
we have
 $$
 \mathcal{H}(\sigma_0(\mu))=\sigma_0(\mu)+(a_2-a_1)\mu=\widehat{\mathcal{G}}(\mu)-a_1\mu.
 $$
arriving at an equivalent definition  of $\widehat{\mathcal{G}}(\mu)$:
\begin{equation}
\label{defi2}
\widehat{\mathcal{G}}(\mu):=\mathcal{H}(\sigma_0(\mu))+a_1\mu.
\end{equation}
Moreover, thanks to the definition of $\sigma_0(\mu)$, applying (\ref{defi1}) we find that
$$
\widehat{\mathcal{G}}(\mu)=\sigma_0(\mu)+a_2\mu\leq \sigma_2+a_2\mu.
$$
Using now, (\ref{defi2}) we have that
$$
\widehat{\mathcal{G}}(\mu)=\mathcal{H}(\sigma_0(\mu))+a_1\mu\leq \sigma_2+a_1\mu.
$$
Hence (\ref{senegal}) follows.

To show that $\widehat{\mathcal{G}}$ is increasing we compute the derivative in (\ref{defi0}) and get
 $$
 -\sigma_0'(\mu)+ \mathcal{H}'(\sigma_0(\mu))\sigma_0'(\mu)=a_2-a_1,
 $$
so that
 $$
 \sigma_0'(\mu)=\frac{a_2-a_1}{ \mathcal{H}'(\sigma_0(\mu))-1}.
 $$
 Therefore,
 $$
 \widehat{\mathcal{G}}'(\mu)=\sigma_0'(\mu)+a_2=\frac{\mathcal{H}'(\sigma_0(\mu))a_2-a_1}{\mathcal{H}'(\sigma_0(\mu))-1}>0,
 $$
 because $\mathcal{H}'(\sigma_0(\mu))<0$.

Finally we prove the limit behaviour of $\widehat{\mathcal{G}}$. First, consider the limit when $\mu\to 0$. It is clear that $\widehat{\mathcal{G}}(0)=\sigma_0(0)=0$.
 Moreover,
\begin{equation}
  \label{eq:limG/mu}
 \lim_{\mu\to 0} \frac{\widehat{\mathcal{G}}(\mu)}{\mu}= \lim_{\mu\to 0} \frac{\sigma_0(\mu)}{\mu}+a_2.
\end{equation}
From~\eqref{defi0} we have
 $$
 \frac{\sigma_0(\mu)}{\mu}\left(\frac{\mathcal{H}(\sigma_0(\mu))}{\sigma_0(\mu)}-1\right)=a_2-a_1.
 $$
 Since $\lim_{s\to 0}\frac{\mathcal{H}(s)}{s}=-\frac{\g_1|\O_2|}{\g_2|\O_1|}$, it follows that
 $$
 \lim_{\mu\to 0} \frac{\sigma_0(\mu)}{\mu}=\g_2|\O_1|\frac{a_1-a_2}{\g_1|\O_2|+\g_2|\O_1|}.
 $$
Hence, replacing this latter expression into~\eqref{eq:limG/mu}
$$
 \lim_{\mu\to 0} \frac{\widehat{\mathcal{G}}(\mu)}{\mu} =\frac{a_1\g_2|\O_1|+a_2\g_1|\O_2|}{\g_1|\O_2|+\g_2|\O_1|}.
$$
On the other hand, to achieve the limit $\widehat{\mathcal{G}}(\mu)\to+\infty$ when $\mu\to\infty$, we will use the behaviour of $\sigma_0(\mu)$ according to whether $a_2>a_1$ or not. Indeed,
if $a_2>a_1$, since $\sigma_0(\mu)\to -\infty$, as $\mu\to +\infty$, using (\ref{defi2}), we conclude that
$$
 \lim_{\mu\to+\infty}\widehat{\mathcal{G}}(\mu)= \lim_{\mu\to+\infty}(a_1\mu +\mathcal{H}(\sigma_0(\mu))= +\infty.
$$
 If $a_2<a_1$ we will use that $\sigma_0(\mu)\to \sigma_2$, as $\mu\to \infty$, and (\ref{defi1}), both proving us with the limiting behaviour of the function $\widehat{\mathcal{G}}(\mu)$ for this particular case.
 Finally, if $a_1=a_2$, the limit follows just using the fact that $\sigma_0(\mu)=0$.

Consequently,
$$
\lim_{\mu\to\infty}\frac{\widehat{\mathcal{G}}(\mu)}{\mu}=\left\{
\begin{array}
  {ll}
 \displaystyle \lim_{\mu\to\infty}\left(a_1+\frac{\mathcal{H}(\sigma_0(\mu))}{\mu}\right)=a_1,\quad& \text{if } a_2>a_1,\\[8pt]
 \displaystyle \lim_{\mu\to\infty}\left(a_2+\frac{\sigma_0(\mu)}{\mu}\right)=a_2,\quad& \text{if } a_2<a_1.
\end{array}
\right.
$$


Now take $\l>\widehat{\mathcal{G}}(\mu)$. If $\l\geq \s_2+a_2\mu$  then
thanks to Lemma~\ref{autovalor}, we first have that
$$ \Lambda_1(-\l+a_1\mu,-\l+a_2\mu)<0.$$
On the other hand, if   $\mu> (\l-\s_2)/a_2$, using again Lemma~\ref{autovalor}, we have that
$$
\begin{array}{rcl}
 \Lambda_1(-\l+a_1\mu,-\l+a_2\mu)<0 & \Longleftrightarrow & \l-a_1\mu>\mathcal{H}(\l-a_2\mu)\\
 & \Longleftrightarrow &(a_2-a_1)\mu>-(\l-a_2\mu) +\mathcal{H}(\l-a_2\mu)\\
  & \Longleftrightarrow &(a_2-a_1)\mu>-\sigma +\mathcal{H}(\sigma),
 \end{array}
$$
where $\sigma=\l-a_2\mu$. Moreover, if $\l>\widehat{\mathcal{G}}(\mu)$ then $\l>\sigma_0(\mu)+a_2\mu$ and hence
$$
\l-a_2\mu>\sigma_0(\mu).
$$
Now, using that  $\mathcal{G}$ is decreasing, we get that $\mathcal{G}(\sigma_0(\mu))>\mathcal{G}(\l-a_2\mu) $, that is, using that $\mathcal{G}(\sigma_0(\mu))=(a_2-a_1)\mu$ we obtain
$$
(a_2-a_1)\mu>-(\l-a_2\mu)+\mathcal{H}(\l-a_2\mu),
$$
then,
$$
\l- a_1\mu>\mathcal{H}(\l-a_2\mu),
$$
or equivalently,   $\Lambda_1(-\l+a_1\mu,-\l+a_2\mu)<0$.
Analogously, we can show that if $\l<\widehat{\mathcal{G}}(\mu)$ then $\Lambda_1(-\l+a_1\mu,-\l+a_2\mu)<0$ and that if $\l=\widehat{\mathcal{G}}(\mu)$ then $\Lambda_1(-\l+a_1\mu,-\l+a_2\mu)=0$.
\end{proof}
As consequence of this result, condition (\ref{condi2igual}) is equivalent to
\begin{equation}
\label{condi2igualequi}
(\mu-g(\l))\cdot(\l-\widehat{\mathcal{G}}(\mu))>0. 
\end{equation}
Finally, we show a non-existence result reciprocal  to Proposition \ref{nomularge}.

\begin{proposition}
\label{nolambdalarge}
Let $\mu>0$ fixed. Then,
there exists $\lambda_*=\lambda_*(\mu)$ such that \eqref{eq:main.system} with $\l_1=\l_2=\l$ does not have coexistence states  if $\lambda>\lambda_*$.
\end{proposition}
\begin{proof} 
Assume by contradiction that there exists a coexistence state $(u_1,u_2, v)$ of \eqref{eq:main.system} for $\lambda$ large enough. Recall that $v\leq \mu$, and then
$$
-\Delta u_i\geq u_i(\lambda -a_i\mu-\alpha_iu_i)\quad\mbox{in $\O_i$.}
$$
Hence,
$$
u_i\geq \vartheta_i\quad \mbox{in $\O_i$,}
$$
where $\vartheta_1$ and $\vartheta_2$ are the positive solutions of (\ref{eq:logisuna}) and (\ref{omega2}) with $\lambda_i=\lambda-a_i\mu$ for $i=1,2$, respectively. 

Consequently,
$$
\mu=\sigma_1^\O[-\Delta+v+b_1u_1\chi_{\Omega_1}+ b_2u_2\chi_{\Omega_2}]\geq \sigma_1^\O[-\Delta+b_1\vartheta_1\chi_{\Omega_1}+ b_2\vartheta_2\chi_{\Omega_2}].
$$
Since,
$$
\vartheta_i\geq\frac{(\lambda-a_i\mu -\sigma_i)}{\alpha_i}\varphi_1^i,\quad i=1,2,
$$
where  $\varphi_1^i$ is the positive eigenfunction associated to $\sigma_i$ with $\|\varphi_1^i\|_\infty=1$  in $\Omega_i$ for $i=1,2$.  Hence,
$$
\mu\geq \sigma_1^\O[-\Delta+b_1\frac{(\lambda-a_1\mu -\sigma_1)}{\alpha_1}\varphi_1^1\chi_{\Omega_1}+ b_2\frac{(\lambda-a_2\mu -\sigma_2)}{\alpha_2}\varphi_1^2\chi_{\Omega_2}].
$$
and passing to the limit as $\l\to\infty$ we arrive at  a contradiction.
\end{proof}

\begin{figure}[ht!]
\centering
\includegraphics[scale=0.25]{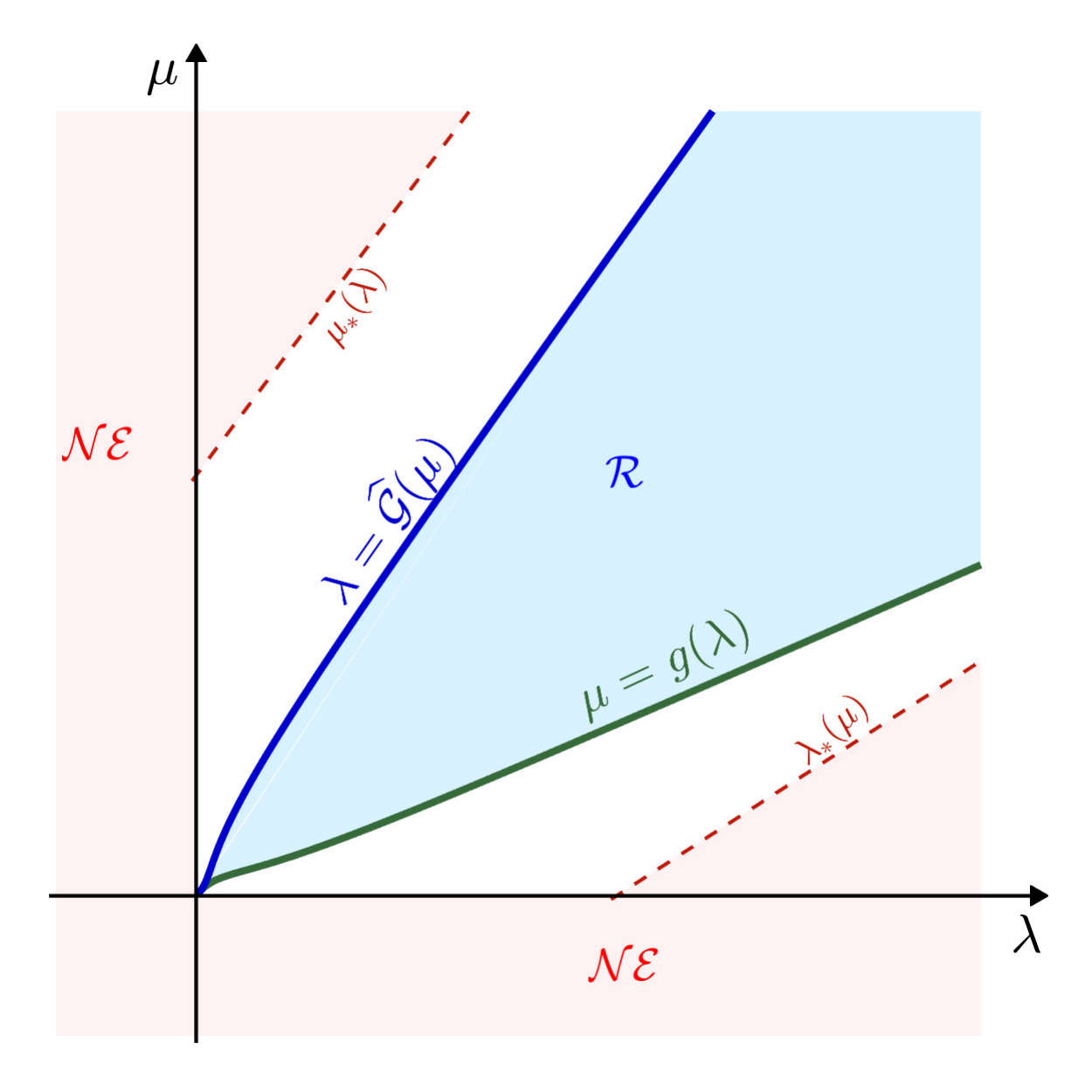}
\includegraphics[scale=0.25]{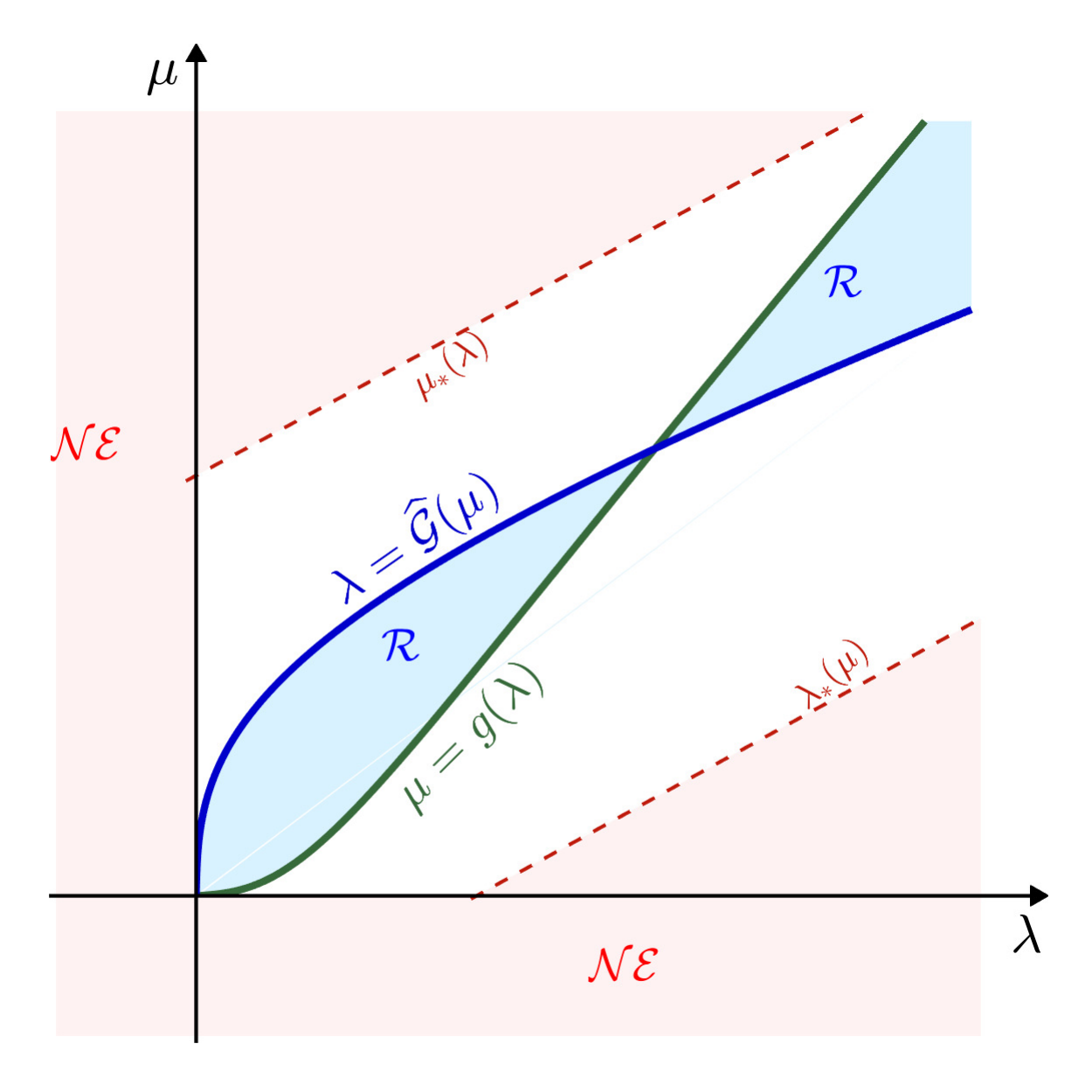}
\caption{We have represented by $\mathcal{R}$ the coexistence region defined by (\ref{condi2igualequi}), that is $\mathcal{R}:=\{(\l,\mu)\in\mathbb{R}^2: (\mu-g(\l))\cdot(\l-\widehat{\mathcal{G}}(\mu))>0\}$. The non-existence region, denoted by $\mathcal{NE}$, includes the results of Corollary \ref{col:nece_cond}, Propositions \ref{nomularge} and \ref{nolambdalarge}. The behaviours of the maps $\mu=g(\lambda)$ and $\lambda=\widehat{\mathcal{G}}(\mu)$ at $0$ and $+\infty$  are given (\ref{limi2}) and (\ref{oro}), respectively}
\label{fig3}
\end{figure}

\section{Limit system}
\label{sect:limit}
Finally we assume that $\mu\geq 0$ and study what will happen in case $\lambda_1$ goes to infinity, for  fixed values of $\lambda_2$ and $\mu$. Let us consider then, $(u_2,v)$, a so-called large-weak solution of the following system
\begin{equation}
\label{eq:large_system2} \left\{\begin{array} {l@{\quad}l}
-\Delta u_{2} =u_{2}(\lambda_2-\alpha_2 u_{2}-a_2 v),& \Omega_2, \\
-\Delta v=v(\mu- v-b_{2}u_{2}),& \Omega_2,\\
u_{2}=\infty,\quad v =0 &\Sigma,\\
\partial_{\bf n} u_{2}=\partial_{\bf n}v= 0,&\partial \Omega,
\end{array}\right.
\end{equation}
in the following sense:  $u_2$ is a classical solution, $u_2\in C^2(\O_2)$, $u_2\geq 0$ in $\O_2$  and the boundary data has to be understood as in~\eqref{eq:large}  and $v\in H^1(\O)$, $v=0 $ in $\O_1$, $v\geq 0$ in $\O_2$ and
$$
\int_{\O_2}\nabla v\cdot\nabla \varphi=\int_{\O_2}v(\mu-v-b_{2}u_{2}) \varphi,\quad\mbox{for all $\varphi\in C^1_c(\Omega_2\cup\Gamma)$,}
$$
defined in (\ref{ho1}.)

\begin{proposition}
\label{Prop:large_limit_1}
Let $\lambda_2$ and $\mu$ be fixed and take a sequence $\l_{1,n}\to +\infty$. Consider a sequence of coexistence states $(u_{1,n},u_{2,n},v_n)$ of \eqref{eq:main.system}  with $\l_1=\l_{1,n}$. Then,  as $\l_{1,n}\to +\infty$
$$
u_{1,n}(x)\to +\infty \quad\mbox{for all $x\in \overline\O_1$,}\quad\quad  u_{2,n}\to u_2 \quad\mbox{in $C^2(\O_2)$},
$$
and
$$
v_n\to v\quad\mbox{weak in $H^1(\Omega_2)$, strong in $L^2(\Omega_2)$ and $C^2(\Omega_2)$.}
$$
\end{proposition}
\begin{proof}
It is clear
that $v_n\leq \mu$. Hence, replacing $v_n$ by $\mu$ in the other two equations of system  \eqref{eq:main.system},
it follows that
$$
-\Delta u_{i,n} \geq u_{i,n}(\lambda_i-\alpha_i u_{i_n}-a_i \mu).
$$
Next, we know, as it is shown in~\cite{AlvarezBrandleMolinaSuarez} that
$$
u_{1,n} \to \infty\quad\mbox{in $\overline\Omega_1$, as $\lambda_1 \to \infty$.}
$$
We study now $v_n$. It is clear that
$$
\int_{\O}|\nabla v_n|^2\leq \mu\int_{\O}v_n^2\leq \mu^3|\O|,
$$
and hence $v_n$ is bounded in $H^1(\O)$, and then $v_n\to v$ weak in $H^1(\O)$ and strong in $L^2(\O)$. Observe that $v_n$ verifies
\begin{equation}
  \label{eq:testvn}
\int_\O\nabla v_n\cdot\nabla\varphi+a_1\int_{\O_1}v_nu_{1,n}\varphi=\int_\O v_n(\mu-v_n)-a_2\int_{\O_2}v_nu_{2,n}\varphi,
\end{equation}
for all $\varphi\in H^1(\Omega)$. We show that $v=0$ in $\O_1$. We argue by contradiction, assume that for $D_1\subset \Omega_1$ we have that $v>0$. Then, taking $\varphi\in C_c^\infty(D_1)$, $\varphi>0$ and we find that
$$
-\int_{D_1} \Delta \varphi v_n+a_1\int_{D_1}v_nu_{1,n}\varphi\leq \mu \int_{D_1} v_n \varphi
$$
such that passing to the limit, and using that $u_{1,n}(x)\to \infty$,  we arrive at a contradiction. Consequently,
$$
v\in H_{0,\Sigma}^1(\Omega_2).
$$
On the other hand, observe that $u_{2,n}\leq L_{\lambda_2}$ in $\Omega_2$, where $L_{\lambda_2}$ is the large solution of problem \eqref{eq:large}. Thus, for any bounded domain $D_2$ inside $\Omega_2$, i.e.
$D_2\subset \Omega_2$ we find that $u_{2,n}$ is bounded in $D_2$. Hence, considering  the equation of $v_n$ in $D_2$,
$$
-\Delta v_n=v_n(\mu-v_n-a_2 u_{2,n})\quad\mbox{in $D_2$,}
$$
we find that $v_n$ is bounded in $C^{2}(D_2)$ and $v_n\to v$ in $C^2(\Omega_2)$.

Furthermore, we consider the second equation
$$
-\Delta u_{2,n}=u_{2,n}(\l_2-\alpha_2u_{2,n}-a_2v_n)\quad \mbox{in $\Omega_2$,}
$$
and then $u_{2,n}\to u_2\in C^2(\Omega_2)$ with
$$
-\Delta u_{2}=u_{2}(\l_2-\alpha_2u_2-a_2v)\quad \mbox{in $\Omega_2$.}
$$
Furthermore, on the interface $\Sigma$ we will arrive at
$$
\frac{\partial u_{2,n}}{\partial {\bf n_2}} +\gamma u_{2,n} =\gamma u_{1,n} \to \infty \quad \hbox{on}\quad \Sigma.
$$
We conclude, see instance \cite{GM-R-S2}, that $u_2=\infty$ in $\Sigma$.


Finally, considering the equation $v_n$ of system  \eqref{eq:main.system}. Thanks to~\eqref{eq:testvn}, taking now a test function $\varphi\in C^1_c(\Omega_2\cup\Gamma)$ it follows that
$$
\int_{\Omega_2}\nabla v_n\cdot\nabla \varphi=\int_{\Omega_2}v_n(\mu-\beta v_n-a_2u_{2,n})\varphi.
$$
Therefore, applying Proposition\;\ref{Prop_v0} we can also say that \eqref{limit:lambda_2} is also satisfied,
just passing to the limit in the third equation of \eqref{eq:main.system} and since $v_n\in C^{2,\eta}(\overline \Omega)$ is positive and bounded in $\overline \Omega$. Hence,
we find the convergence of $v_n$ to $v$ in $\Omega_2$ while
vanishing in $\Omega\setminus\Omega_2$.
\end{proof}




\end{document}